\documentclass[times,sort&compress,3p]{elsarticle}
\journal{Journal of Multivariate Analysis}
\usepackage[labelfont=bf]{caption}

\usepackage{amsmath,amsfonts,amssymb,amsthm,booktabs,color,epsfig,graphicx}
\PassOptionsToPackage{hyphens}{url}\usepackage{hyperref}
\usepackage{mathrsfs}
\DeclareMathAlphabet{\mathpzc}{OT1}{pzc}{m}{it} 
\usepackage{mathbbol}
\usepackage{calligra}

\usepackage[T1]{fontenc} 

\theoremstyle{plain}
\newtheorem{theorem}{Theorem}

\theoremstyle{definition}
\newtheorem{definition}{Definition}
\newtheorem{remark}{Remark}
\newtheorem{example}{Example}

\usepackage{pgfplots}
\pgfplotsset{compat = newest}

\usepackage{wrapfig}

\usepackage{cancel}

\newcommand{\Rbb}{\mathbb{R}}

\newcommand{\Kbb}{\mathbb{K}}
\newcommand{\Prbb}{\mathscr{P}}

\newcommand{\Hcal}{\mathcal{H}}

\newcommand{\Scal}{\mathcal{S}}

\newcommand{\Cscr}{\mathscr{C}}

\DeclareMathOperator*{\argmax}{arg\,max}
\makeatletter
\newcommand*{\defeq}{\mathrel{\rlap{%
                     \raisebox{0.3ex}{$\m@th\cdot$}}%
                     \raisebox{-0.3ex}{$\m@th\cdot$}}%
                     =}
\makeatother

\newcommand{\say}[1]{``#1''}

\begin{document}

\begin{frontmatter}

\title{Data depth functions for non-standard data by use of formal concept analysis}

\author[1]{Hannah Blocher\corref{mycorrespondingauthor}}
\author[1]{Georg Schollmeyer}

\address[1]{Department of Statistics,
Ludwig-Maximilians-Universität München,
Ludwigstr. 33,
80539 München,
Germany}

\cortext[mycorrespondingauthor]{Corresponding author. Email address: \url{hannah.blocher@stat.uni-muenchen.de}}

\begin{abstract}

In this article we introduce a notion of depth functions for data types that are not given in standard statistical data formats. We focus on data that cannot be represented by one specific data structure, such as normed vector spaces. This covers a wide range of different data types, which we refer to as non-standard data. 
Depth functions have been studied intensively for normed vector spaces. However, a discussion of depth functions for non-standard data is lacking. In this article, we address this gap by using formal concept analysis to obtain a unified data representation. Building on this representation, we then define depth functions for non-standard data. Furthermore, we provide a systematic basis by introducing structural properties using the data representation provided by formal concept analysis. Finally, we embed the generalised Tukey depth into our concept of data depth and analyse it using the introduced structural properties. Thus, this article presents the mathematical formalisation of centrality and outlyingness for non-standard data and increases the number of spaces in which centrality can be discussed. In particular, we provide a basis for defining further depth functions and statistical inference methods for non-standard data. 
\end{abstract}


\end{frontmatter}

\section{Introduction}\label{sec:intro}

\textit{Data depth functions} generalise the concept of \textit{centrality} and \textit{outlyingness} to multivariate data and provide therefore a useful concept to define \textit{nonparametric} and \textit{robust statistical methods}. To achieve this, depth functions denote the center and outlying areas based on an underlying distribution or a data cloud. Moreover, they classify how near other areas are to the outlying or center ones. Since there does not exist one unique perspective on centrality, several different depth functions on $\mathbb{R}^d$ have been developed. Some examples are simplicial depth, see \cite{liu90}, zonoid depth, see \cite{dyckerhoff96}, or Tukey depth, see \cite{tukey75, donoho92}. To give a systematic basis on the notion of depth functions \cite{zuo00} and \cite{mosler13} set up properties to mathematically formalise the existing intuition about centrality and outlyingness in $\mathbb{R}^d$. In addition to $\mathbb{R}^d$, the scope of data types on which depth functions are defined increased in the last years. For example, in \cite{gijbels17} the authors discuss properties of depth functions on functional data. \cite{goibert22} developed a depth function for total orders. \cite{geenens23} went the next step and defined depth functions for general metric spaces.
In recent years, this concept has been used to construct statistical methods and to analyse different application settings. For example, \cite{li04, pawar22} developed statistical nonparametric and robust tests using depth functions. Using the robust central-outward order \cite{mozharovskyi22} studied anomaly detection. In \cite{chebana11}, data depth functions are used to visualise and describe features of hydrologic events given two real world streamflow data sets from Canada.

However, all these depth functions and analyses have in common that they build on spaces that have a strong underlying structure, like Banach spaces. Thus, to apply the concept of depth and resulting statistical methods, the data has to be embedded into statistical standard data formats, like numeric or spatial. 
The aim of this article is to consider depth functions without assuming one specific data type in advance. This includes data types like, e.g., set of partial orders or mixed spatial and ordinal data. We call data that cannot be embedded into statistical standard data formats \textit{non-standard data}. To give a notion of centrality and outlyingness for non-standard data, a unified, flexible and general applicable data representation is needed. Therefore, we use the theory of \textit{formal concept analysis} which transforms the data set into a \textit{closure system} on the entire data set itself. In particular, it does not force the user to make assumptions about the data that are not necessarily true, simply for the sake of mathematical feasibility. We use this unified representation of the data to generally define depth functions for non-standard data. To clarify the notion of centrality and outlyingness for data represented via formal concept analysis we introduce \textit{structural properties}. Thereby, we use the structure of the closure system to transfer the properties given by \cite{zuo00} and \cite{mosler13} and develop further properties representing the connections between the data points. In this context, we investigate that, unlike $\mathbb{R}^d$, non-standard data can have data elements that are naturally outlying or central based on the structure of the underlying space. Finally, we embed the \textit{generalised Tukey depth} given in \cite{blocher22} into our concept of depth functions for non-standard data. Moreover, we use the provided structural properties to analyse the generalised Tukey depth. Note that, analogous to the generalised Tukey depth, the concrete outlier measures developed in~\cite{hu23} can also be analysed using the structural properties introduced here. These outlier measures are also based on formal concept analysis. Therefore, only a slight adaptation to the concept of depth introduced here is necessary.

This article is intended to give a systematic basis on depth functions for non-standard data. By introducing structural properties we give a notion of centrality and outlyingness. Furthermore, we provide a framework to analyse depth functions for non-standard data and start a discussion on the notion of centrality for such kind of data.

The structure of this article is as follows: We begin with concrete motivation examples. Since formal concept analysis is the basis of our general notion of a depth function, we give a short introduction to the theory of formal concept analysis in Section~\ref{sec:preliminaries}. Based on this, Section~\ref{sec:definition} introduces the general definition of depth functions for non-standard data by use of formal concept analysis. There we define the Tukey depth within our concept of depth functions. The next section states the structural properties and their restriction strength. Afterwards, we analyse the generalised Tukey depth using structural properties. In Section~\ref{sec:conclusion}, we collect our concluding remarks.

\section{Motivation}\label{sec:motivation}

The purpose of this section is to present motivating examples to show that depth functions based on formal concept analysis cover a wider range of data types than the currently discussed depth functions known to the authors. Instead of using observations directly to define a depth function, we first apply formal concept analysis. This gives us a definition of depth functions that can be applied to a wide variety of data types simultaneously. Most importantly, it is simple for the practitioner to use, as it only relys on the definition of a cross-table, as will be seen in this section.

One of the goals of formal concept analysis is to discover relational connections between observations. Therefore, the observations are grouped according to the observed values. For example, consider the Titanic data set, where each passenger's age, sex and the passenger's class of travel is recorded. One group might be all passengers aged between 20 and 40 who are female. Another smaller group are all female passengers between the ages of 20 and 40 travelling in first class. In this case, the smaller group is more specific, e.g. in the smaller group the passengers also fulfil the condition that they all travel in first class. Formal concept analysis formalises the arrangement into different groups. The first step is to define a formal context that represents the data set. In a nutshell, a formal context is a generalisation of a cross-table. 

To provide an intuitive introduction to the main concepts of formal concept analysis and to illustrate the motivation for defining depth functions based on formal concept analysis, we give now two examples. These examples also demonstrate the variety of different data types that can be represented by a formal context.

\begin{example}\label{motivation_titanic}
For the Titanic data example provided by the Titanic machine learning competition, see~\cite{cukierski12}, we get the formal context in Table~\ref{tab:titanic_motivation} (this is just a snippet). Here each row represents an observation/passenger (i.e. called \textit{objects} in formal concept analysis). The columns are binary \textit{attributes} that can be either true or false. Transforming the non-binary observations (such as age, spatial values, etc.) into binary attributes is called \textit{conceptual scaling}, see~\cite[p. 36ff]{ganter12}. Since the number of attributes can be infinite, this is a flexible, but also accessible approach to data representation. For classical data types there are already standard scaling methods, see~\cite[p. 36ff]{ganter12}. As an example, consider  Table~\ref{tab:titanic_motivation}. The first two columns represent the sex of the passenger (f: female and m: male). The next three columns are the class level of the passenger (I, II and III). These two nominal observations are included via \textit{nominal scaling}, see~\cite[p. 42]{ganter12}. The remaining columns represent the age of the passenger, e.g. $g_1$ is 34.5 years old. The corresponding attributes are the statements \say{$\le x$} and \say{$\ge x$} for all $x \in \mathbb{R}$ with $(g, \text{\say{$\le x$}}) \in I$ if and only if the numeric observation of $g$ is smaller or equal to $x$. Analogously \say{$\ge$} is defined. This scaling method is called \textit{interordinal scaling}, see~\cite[p. 42]{ganter12}. Note that this scaling method implies that there exists no object that has every attribute. The crosses in Table~\ref{tab:titanic_motivation} indicate that the attribute applies to the passenger/observation. How different observation values are represented as a set of binary attributes depends on the scaling method used, see~\cite[p. 36ff]{ganter12}.

We obtain the groups that represent the relationship between the observations by considering all possible combinations of attributes and summarising all the observations that agree on all these attributes. Note that these groups may overlap. They also have a certain mathematical structure, i.e. they show the implications between the passengers/observations. For example, in Table~\ref{tab:titanic_motivation} we see that if $g_3$ and $g_4$ are in a group, we can directly imply that $g_2$ is also in the group, since there is no group with $g_3$ and $g_4$ that does not also contain $g_2$. By collecting all these implications, we can fully determine the arrangement of the groups given by the formal context. For further details see Section~\ref{sec:preliminaries}.
\begin{table}[!ht]
\begin{tabular}{|l||c|c|c|c|c|c|c|c|c|c|c|c|c|c|l|}
\hline
       & f & m & I & II & III & $\ldots$ & `$\le 23$' & $\ldots$ & `$\le 24$' & $\ldots$ & `$\le 34.5$' & $\ldots$ & `$\le 35$' & $\ldots$ & {`$\le 47$'} \\ \hline \hline
$g_1$  &   & x &   &    & x   &          &            &          & 	        & 		   & x            & $\ldots$ & x          & $\ldots$ &  x           \\ \hline
$g_2$ & x &   &   &    & x   &          &            &          &            &          &              &          &            &          &  x           \\ \hline
$g_3$  & x &   &   & x  &     &          &            &          &            &          &              &          &            &          &               \\ \hline
$g_4$ & x &   & x &    &     &          & x          & $\ldots$ & x          & $\ldots$ & x            & $\ldots$ & x          & $\ldots$ &  x            \\ \hline
$g_5$  &   & x &   & x  &     &          &            &          &            &          &              &          & x          & $\ldots$ &  x            \\ \hline
$\ldots$  & $\ldots$  & $\ldots$ & $\ldots$  & $\ldots$  &  $\ldots$   & $\ldots$ & $\ldots$ & $\ldots$ & $\ldots$ & $\ldots$ & $\ldots$ &  $\ldots$  & $\ldots$ & $\ldots$ & $\ldots$ \\ \hline
\end{tabular}
\vspace{1em}

\begin{tabular}{|l||c|c|c|c|c|c|c|c|c|c|c|c|c|c|}
\hline
       & $\ldots$ & `$\le 67$' & $\ldots$ & $\ldots$ & `$\ge 67$' & $\ldots$ & `$\ge 47$' & $\ldots$ & `$\ge 35$' & $\ldots$ & `$\ge 34.5$' & $\ldots$ & `$\ge 24$' & $\ldots$ \\ \hline\hline
$g_1$  & x        & x        & x        &          &            &          &            &          &            &          & x            & x        & x          & x        \\ \hline
$g_2$ & x        & x        & x        &          &            &          & x          & x        & x          & x        & x            & x        & x          & x        \\ \hline
$g_3$  &          & x        & x        &          & x          & x        & x          & x        & x          & x        & x            & x        & x          & x        \\ \hline
$g_4$ & x        & x        & x        &          &            &          &            &          &            &          &              &          &           &         \\ \hline
$g_5$  & x        & x        & x        &          &            &          &            &          & x          & x        & x            & x        & x          & x        \\ \hline
$\ldots$  & $\ldots$  & $\ldots$ & $\ldots$ & $\ldots$ & $\ldots$ & $\ldots$ & $\ldots$ &  $\ldots$ & $\ldots$ & $\ldots$ & $\ldots$ & $\ldots$ & $\ldots$ & $\ldots$ \\ \hline
\end{tabular}

\caption{This formal context contains five observations/passengers of the Titanic. Each row represents a passenger. 
The data is from the Titanic machine learning competition, see~\cite{cukierski12}, with PassengerIds 892, 893, 894, 904 and 908. Note that the age $34.5$ is due to the encoding notation. With~\say{$\ldots$}, we denote that the attribute set is indeed infinite and not all attributes are denoted within the table. The same is true for the objects as this is only a snippet of the titanic dataset.}
\label{tab:titanic_motivation}
\end{table}

We emphasise that the choice of binary attributes, as well as the decision which observation has an attribute and which does not, results in a concrete, easily accessible representation of the data. In particular, it exposes many implicit assumptions about the data. Moreover, the data are represented in such a way that we do not impose assumptions that are not present in the data, e.g. that the data can be embedded in a normed vector space. However, we want to point out that the choice of attributes can have a strong influence on the arrangement of the groups and therefore on the depth values. Instead of using interordinal scaling for the age component in the Titanic data set, another approach called ordinal scaling, which includes only the upper bound `$\le x$', can be used, see~\cite[p. 42]{ganter12}. With this smaller formal context, age can only influence the groups through its upper bound. Thus, instead of looking at all passengers between the ages of, e.g., 20 and 40, with this smaller formal context, only all passengers with an age smaller than or equal to 40 are grouped without a lower bound.

We now use the structure given by the different groups and the resulting implications to give a notion of depth to the observations. In the Titanic example, we want to rank the five passengers according to how central or outlying they are compared to the other passengers. In the spirit of \cite{serfling00}, in this paper we introduce properties that represent the notion of centrality and outlyingness.  For example, Property~(P1) states that if we use different scaling methods of the observations (i.e., use different binary attributes in the formal context) that lead to the same groups, the depth function should stay the same. As another example, consider the relationship between $g_2$, $g_3$, and $g_4$. Since any group that includes $g_3$ and $g_4$ also includes $g_2$, we have that $g_2$ has all the properties that $g_3$ and $g_4$ share. So we can say that $g_2$ is at least as specific as $g_3$ and $g_4$ combined. Following the idea of the quasiconcavity property in $\mathbb{R}^d$, see~\cite{mosler13}, Property~(P7i and ii) states that the depth of $g_2$ must be at least as high as the minimum depth of $g_3$ or $g_4$. Other properties discuss that non-standard data may inherit a natural center-outward order, which is not present in $\mathbb{R}^d$. For example, Property~(P4) states that an object lying in every group must be central, regardless of the probability measure.
\end{example}

\begin{example}\label{motivation_posets}
Another example that illustrates the wide range of data types that can be discussed using formal concept analysis is the set of partial orders. A data example for partial orders is given in~\cite{blocher23}, where the objects are partial orders of classifiers. In that article, we considered a set of classifiers applied to a set of data sets. These classifiers were evaluated for each individual data set by several performance measures. Then an algorithm $i$ outperforms another algorithm $j$ on a given data set if and only if there exists one performance measure that states that algorithm $i$ is better than algorithm $j$ and all other performance measures agree that algorithm $i$ is not worse than algorithm $j$. This gives us a partial order for each single data set that represents the performance structure of the classifiers for that data set. Thus, a benchmark suite of 80 data sets results in 80 observed partial orders. Now, to represent this kind of data, in~\cite{blocher22} we introduced a formal context for partial orders in general. Here, each row corresponds to a partial order. The columns of the formal context are the binary attributes needed to describe the partial order. In~\cite{blocher22}, the first $n\cdot(n-1)$, where $n$ is the number of classifiers/items being compared, denotes all possible pairs where we have stated that one classifier dominates the other. The next $n\cdot(n-1)$ rows indicate the reverse, i.e. that one classifier does not dominate another. Note that including the second part as binary attributes is by no means straightforward. But as in~\cite{blocher22, blocher23} we decided to state that non-existence of dominance is indeed a precise observation and we do not want a push towards the linear extensions of the partial orders (which exists when deleting the second $n\cdot (n-1)$ columns), we have therefore included this part. 
\end{example}

As can be seen in the two examples above, and also in Example~\ref{hierarchical nominal} and~\ref{exampl: spatial_1} below, defining a formal context gives us a flexible tool for evaluating data ranging from mixed nominal and ordinal observations to spatial data or partial orders. These examples demonstrate that formal concept analysis clearly reveals the relationships between data points, providing a simple and unified view of the underlying data structure. Finally, we want to highlight differences and connections between our approach and other existing approaches. In contrast to \cite{goibert22, geenens23}, we have completely turned away from a metric approach in our definition of depth functions. \cite{goibert22} discusses total orders and therefore a translation to our approach is straightforward, just restrict the space of all partial orders to the total orders, see Example~\ref{motivation_posets}. This is not so clear for~\cite{geenens23} as the representation via a formal context depends on the concrete structure of the space. Example~\ref{exampl: spatial_1} gives an approach for $\mathbb{R}^d$ equipped with the Euclidean distance. 

\section{Formal concept analysis}\label{sec:preliminaries}

This section briefly describes the parts of formal concept analysis (FCA) required for this article and illustrates them using the two examples above. It is based on \cite{ganter12}. For further readings, we refer to \cite{carpineto04}. FCA was developed by Rudolf Wille, Bernhard Ganter and Peter Burmeister to build a bridge between mathematical lattice theory and applied users. It enables the analysis of relationships between the data points by representing the data in a unified and user-friendly manner.

The fundamental definition of FCA is the representation of a data set as a cross-table, see \cite[p. 17]{ganter12}.
\begin{definition}
	A \textit{formal context} $\Kbb$ is a triple $(G,M,I)$ with $G$ being the \textit{object} set, $M$ the set of \textit{attributes} and $I \subseteq G\times M$ a binary relation between $G$ and $M$. 
\end{definition}
In our case, the objects are the data points and the attributes are characteristics of these data points. The relation $I$ then states whether an object $g$ has an attribute $m$, if $(g,m) \in I$, or not, if $(g,m) \not\in I$.  Thus, these attributes need to be binary-valued, whether they occur or not. Naturally, there exist characteristics of the data points which are many-valued, like sex or age. To include these many-valued characteristics as well into the formal context, we use \textit{conceptual scaling methods}, see \cite[p. 36ff]{ganter12} and Section~\ref{sec:motivation}, which transfers many-valued characteristics into a set of binary-valued attributes. Examples of a scaling methods are given in the Titanic snippet in Example~\ref{motivation_titanic}.
%
By using scaling methods, we can represent a large variety of different data sets through a formal context. This allows us to transfer the most diverse data types into a uniform structure. Also, data sets which are not given in standard statistical data formats. We call such data \textit{non-standard data}. Examples are given in Example~\ref{motivation_titanic}, \ref{motivation_posets}, \ref{hierarchical nominal} and~\ref{exampl: spatial_1}. 
The scaling method can also be used to reduce data complexity. This can lead to a conceptual scaling error, see for example \cite{hanika21}.

Based on this user-friendly representation of the data set by a formal context, we can now define so-called \textit{derivation operators}, see \cite[p. 18]{ganter12}:
\begin{align*}
	&\Psi: 2^G \to 2^M, A \mapsto  A' \defeq \{m \in M \mid \forall g \in A \colon (g,m) \in I\},\\
	&\Phi: 2^M \to 2^G, B \mapsto B' \defeq \{g \in G \mid \forall m \in B \colon (g,m) \in I\}.
\end{align*}
The function $\Psi$ maps a set of objects $A$ onto every attribute which every object in $A$ has. So, $A' = \Psi(A) = \cap \{m \in M \mid (g,m) \in I\}$ is the maximal set of attributes that every object in $A$ has. The reverse, from attribute set to object set, is provided by the function $\Phi$. We set $\Phi(\emptyset) = G$ and $\Psi(\emptyset) = M$. The composition of these two functions $\gamma \defeq\Phi \circ \Psi: 2^G \to 2^G$ gives us then a family of sets which denotes the relationship between the data points. We have $\gamma(A) = \Phi \circ \Psi(A) = \cap_{m \in \Psi(A)}\{g \in G \mid (g,m) \in I\}$. In other words, every object set $E \subseteq G$ which is an element of the codomain of $\gamma = \Phi \circ \Psi$ is the maximal set of objects which have all the same attributes $\Psi(E)$ in common. Thus, the composition groups all those objects together that have the same attributes. With slight abuse of notation, for any $m \in M$  we write $\Phi(m)$ for $\Phi(\{m\})$. The same convention holds for $\Psi(g)$ and $\gamma(g)$ with $g \in G$.

\begin{definition}
	The set $\gamma(2^G)$ is called the set of \textit{extents}, see \cite[p. 18]{ganter12}. The set $\gamma(A)$ is mentioned as \textit{extent set}, shortly \textit{extent}, of $A\subseteq G$.
\end{definition}

To take advantage of this slightly different representation of the data set as a family of sets, we use that $\gamma$ defines a \textit{closure operator}, see \cite[p. 8]{ganter12}. In what follows, we use the term \textit{closure} always in the context of a closure operator or a closure system. When referring to a closed set based on a topology, metric, or norm we denote this by \textit{topological(-ly) closed/closure}.
\begin{definition}
	A closure operator $\gamma: 2^G \to 2^G$ is defined as a function on a power set to itself. A closure operator needs to be extensive (i.e. for all $A \subseteq G, \: A \subseteq \gamma(A)$), isotone (i.e. if $A\subseteq B \subseteq G$, then $\gamma(A) \subseteq \gamma(B)$) and idempotent (i.e. for all $A \subseteq G, \: \gamma(A) = \gamma(\gamma(A))$). 
	
	In particular, a closure operator always induces a \textit{closure system} $\gamma(2^G)$. 
	A closure system $\Scal \subseteq 2^G$ is a family of sets which contains the entire space (i.e. $G \in \Scal$) and any intersection of sets in $\Scal$ is again in $\Scal$ (for all $S \subseteq \Scal$ with $S \neq \emptyset$ we have $\bigcap_{s \in S} s \in \Scal$).
\end{definition}
Note that there exists a one-to-one correspondence between the closure system and the closure operator. Since $\gamma = \Phi \circ \Psi$ is a closure operator, the set of extents is a closure system. Thus, the closure operator $\gamma$ describes the closure system and vice versa. For more details on closure systems see \cite[Chapter 0.3]{ganter12}. 

\begin{example}\label{exampl: nominal_ordinal_2}
	Recall the Titanic Exampe~\ref{motivation_titanic}. We obtain $\Psi(g_2)= \{\text{f}, \text{III}, \text{`}\le 47\text{'}, \ldots, \text{`}\le 67\text{'}, \ldots, \text{`}\ge 47\text{'}, \ldots, \text{`}\ge 35\text{'} \}.$ Therefore, we get $\gamma(g_2) = \{g_2\}$. One can show that the set of extents $\gamma\left(2^{\{g_1, g_2, g_3, g_4, g_5\}}\right)$ equals \begin{align*}
	&\{\emptyset, \{g_1\}, \{g_2\}, \{g_3\}, \{g_4\}, \{g_5\},
	\{g_1, g_4\}, \{g_1, g_5\}, \{g_2, g_5\}, \{g_2, g_3\}\\
	\{g_1, g_4, g_5\}, & \{g_1, g_2, g_5\}, \{g_2, g_3, g_5\},
	\{g_1, g_2, g_4, g_5\}, \{g_1, g_2, g_3, g_5\},
	\{g_1, g_2, g_3, g_4, g_5\}\}.	
	\end{align*}
\end{example}

Let us take a closer look at how the closure operator describes the connection between the data points. This is the basic idea of how we will later use the closure operator to define structural properties for depth functions. As we pointed out above the closure operator groups data points together which have certain attributes in common. This means if $a \in \gamma(A) \setminus A$ the object $a$ has all attributes which every object in $A$ has as well. Thus, one can say that $A$ \textit{implies} $a$ based on the relationship structure given by the formal context which is then included in the definition of $\gamma$. Therefore using the closure system or closure operator (both describe the same since they have a one-to-one correspondence) to define the structural properties of the depth function illustrates the relationship between the data points. For further details on implications see \cite[Chapter 2.3]{ganter12}. Note that in \cite{ganter12} attribute implications are discussed and we focus here on object implications. Nevertheless, the concepts can be transferred to object implications.


\begin{example}\label{hierarchical nominal}

	Now, we introduce a further scaling method, the so-called \textit{hierarchical nominal scaling}. This scaling method is inspired by the occupations of persons within a social survey. 
	Usually occupations are categorised within a hierarchy of different levels, for example within the International Standard Classification of Occupations (ISCO) of 2008, see \url{https://www.ilo.org/publications/international-standard-classification-occupations-2008-isco-08-structure} (accessed: 29.07.2024). On a first level occupations are split into different categories $(a_1,b_1,...)$. For example category $a_1$ could be \say{Managers}, category $b_1$ \say{Professionals}, etc. Each of these categories is then split again on a more fine-grained level (Level~2) into further subcategories. In this case, category $a_1$ is split into $a_1a_2$ \say{Managers: Chief executives, senior officials, and legislators} and $a_1b_2$ \say{Managers: Administrative and commercial managers} and so on. Note that the Level~2 subcategories based on the first level $b_1$ split do not have to match those of the Level~2 splits based on the Level 1 $a_1$ split. Subsequently, the Level~2 categories are again subdivided into subcategories, and so on. Such data structure can be conceptually scaled in a natural way. For every level, we introduce attributes describing every single category based on the upper level classification, see Table~\ref{crosstable_hierarchical_nominal}. Here, for simplicity, we used only two levels with two categories, respectively.

	\begin{table}[ht]
		\centering
		\begin{tabular}{r|rr|rrrr}
			\hline
			& $a_1$ & $b_1$ & $a_1a_2$ & $a_1b_2$ & $b_1a_2$ & $b_1b_2$ \\ 
			\hline
			$a_1a_2$ & x &  & x &  &  &  \\ 
			$a_1b_2$ & x &  &  & x &  &  \\ 
			$b_1a_2$ &  & x &  &  & x &  \\ 
			$b_1b_2$  &  & x &  &  &  & x \\ 
			\hline
		\end{tabular}
		\caption{\label{crosstable_hierarchical_nominal} Minimal example of a hierarchical nominal scaling with two levels ($1$, $2$) and two categories ($a$,$b$, respectively).}
	\end{table}
Now, let us take a closer look at the extents given by Table~\ref{crosstable_hierarchical_nominal} which are $$\{\emptyset, \{a_1a_2\}, \{a_1b_2\}, \{b_1a_2\},\{b_1b_2\}, \{a_1a_2, a_1b_2\},\{b_1a_2, b_1b_2\}, G\}.$$ Furthermore, we obtain that $a_1a_2 \in \gamma(\{a_1b_2, b_1a_2\}) = \Phi \circ \Psi (\{a_1b_2, b_1a_2\}) = \Phi(\emptyset) = G$.
\end{example}

\begin{example}\label{exampl: spatial_1}
	In the next sections, we want to transfer the idea of depth functions from $\mathbb{R}^d$ to general non-standard data which are represented by a formal context. 
	Before looking at this, we now show how one can represent the elements in $\mathbb{R}^d$ as objects of a formal context. We consider $\mathbb{R}^d$ together with the topology induced by the Euclidean norm. The scaling method introduced here is inspired by \cite{schollmeyer17a, schollmeyer17b}. Let $G = \mathbb{R}^d$ be the object set and the attribute set $M = \{H \subseteq \mathbb{R}^d \mid H \text{ topologically closed halfspace}\}$. Define the relation $I$ between $M$ and $G$ by $(g,H) \in I \Leftrightarrow g \in H$, and let $\Kbb = (G,M,I)$ denote the corresponding formal context. 
	
	With this definition, let's consider $\gamma(2^G)$ induced by $\Kbb$. Let $A \subseteq G = \mathbb{R}^d$, then $\Psi(A)$ are all halfspaces containing every object/point in $A$. Further on, $\gamma(A) = \Phi\circ\Psi(A)$ are then every object/point which lies in every halfspace in $\Psi(A)$. Thus $\gamma(A)$ is the intersection of all halfspaces in $\Psi(A)$ and therefore a topologically closed convex set. 
	More generally, one can show that for every topologically closed convex set in $2^{\mathbb{R}^d}$ there exists a set $A\subseteq G =\mathbb{R}^d$ such that this convex set equals $\gamma(A)$. Thus, $\gamma$ is the convex closure operator on $\mathbb{R}^d$ and the extent set $\gamma(2^G)$ is the set of all topologically closed convex sets.
	
	More concretely, assume that $d = 2$ and $f \in \gamma(\{a,b,c\})$. This means that $f$ lies in every topologically closed halfspace which contains also $a,\: b$ and $c$. In other words, $f$ shares the same attributes as $a,b,c$ share. Thus, one can say that $f$ is implied by $a,b$ and $c$ based on $\gamma$. If $e \not\in \gamma(\{a,b,c\})$, there exists a halfspace which contains $a,b$ and $c$ but not $e$. To summarise, an object $g$ lies in $\gamma(\{a,b,c\})$ if and only if $g$ lies within the triangle given by the vertices $a$, $b$ and $c$. This shows how the concrete definition of the closure system enhances the connection between single objects.
\end{example}

Finally, as we have now introduced several different data situations with different scaling methods, see Example~\ref{motivation_titanic} (nominal, numeric data), Example~\ref{motivation_posets} (partial orders), Example~\ref{hierarchical nominal} (hierarchical nominal data) and Example~\ref{exampl: spatial_1} (spatial data), we want to point out the importance of the scaling method used. First of all, deciding which observations (age, spatial, etc.) to include in the analysis is a general issue in statistics. At first sight, FCA may seem to add a further difficulty by requiring the definition of binary attributes. However, FCA provides a clear understanding of the relationships between observations. Thus, to check whether the scaling method represents the objective, one can look at all the extents as well as all the implications. If the groups do not represent an important fact about the data (e.g. the non-dominance part of the partial orders or the lower bound of the age value), then more attributes are needed. Conversely, if the extents are close to the power set of $G$, one should examine the data set to see if there really are no implications between the data points. If so, some sets of attributes need to be replaced or even deleted. While nominal, ordinal and interordinal scaling are commonly used scaling methods, it is important to discuss carefully for each data set how fine/coarse the arrangement of the groups should be and to consider several different scaling methods.

\section{Definition of depth functions for non-standard data using formal concept analysis}\label{sec:definition}
Our aim in this section is to give a general definition of data depth functions for non-standard data using FCA. By representing the data points $G$ via a formal context, with $G$ being the object set, we obtain a unified structure that is not tailored to one specific data type. In particular, the newly provided depth function allows to analyse a large variety of different data types on centrality and outlyingness issues. With this, nonparametric methods can be developed for all these data.

The depth function presented here only specifies the domain and codomain but not the exact mapping rule. Thus, the structural properties presented later, see Section~\ref{sec:des_properties}, can be seen as generic properties for this kind of depth functions. 

\begin{definition}\label{def:allg_depth_vorschrift}
We define a depth function using FCA by
\begin{align*}
	D_G:G \times \varkappa_G \times \Prbb_G\to \mathbb{R}_{\ge 0}
\end{align*}
for a fixed set of objects $G$ and a set of formal contexts $\varkappa_G  \subseteq \{\mathbb{K} \mid G \text{ is object set of } \mathbb{K}\}$. $\Prbb_G$ is a set of probability measures on $G$ defined on a $\sigma$-field which contains all extent sets of formal contexts in $\varkappa_G$.
\end{definition}

Thus, we compute the depth of an object set based on a probability measure and formal context representing the object relationships. 
We want to emphasise that $G$, $\varkappa_G$ and $\Prbb_G$ depend on each other. Note that Definition~\ref{def:allg_depth_vorschrift} allows the restrictions of $\varkappa_G$ and $\Prbb_G$ to a subset of all possible formal contexts or probability measures, which is sometimes necessary. 
For example, assume that $G = [0,1]$, then for every subset of $G$ there exists a formal context such that this subset is an extent. Thus, in this case, a restriction to a subset of contexts is necessary if we want to allow the uniform distribution to be an element of $\Prbb_G$. This follows from \cite[Chapter 1.1]{tao11} which shows that there cannot exist a probability measure on $[0,1]$ which formalises the intuition of volume and has the entire power set of $[0,1]$ as input. Another aspect is that restriction can lead to a set of formal contexts fulfilling additional structural requirements. With this, it can be possible to define mapping rules which are not possible in general. (See Section~\ref{sec:des_properties}: Property (P8) does not hold for every formal context, but only for a subset.) Another example is given in \cite{blocher22} and \cite{blocher23} where we considered one single formal context on the set of partial orders, see also Example~\ref{motivation_posets}. There, we used the structure given by this concrete formal context to define the mapping rule. 
The same reasoning can be applied to the probability set $\Prbb_G$.
Thus, for a proper definition of the depth function, not only the exact mapping rule is important, but the considered formal contexts and probability measures as well. 

A data set can be represented by different attribute sets and corresponding binary relations. Thus, Definition~\ref{def:allg_depth_vorschrift} can lead to many different depth functions even if the object set $G$, the probability $\Pr \in \Prbb_G$ and a concrete mapping rule are specified. Hence, the choice of formal context and scaling method can have a huge impact on the depth values, i.e. this can also be seen in Section~\ref{sec:des_properties}.

The empirical depth function corresponds to the depth function in Definition~\ref{def:allg_depth_vorschrift} with the empirical probability measure as input. To ensure that the empirical depth function is well defined we need to assume that every empirical probability measure $\Pr_G^{(n)}$ of every probability measure $\Pr_G \in \Prbb_G$ is again an element of $\Prbb_G$. 

\begin{definition}\label{def:allg_empir_depth_vorschrift}
Let $G$ be a set and $\varkappa_G  \subseteq \{\mathbb{K} \mid G \text{ is object set of } \mathbb{K}\}$ a set of formal contexts on $G$. We assume that $\Prbb_G$ consists of probability measures that are defined on a $\sigma$-field containing all extents of $\varkappa_G$. Furthermore, for every sample $g_1, \ldots, g_n$ based on a probability measure $\Pr_G \in \Prbb_G$, we have $\text{Pr}^{(n)}_G \in \Prbb_G$, where $\text{Pr}^{(n)}_G $ is the empirical probability measure based on the sample. Then the empirical depth function for a sample $g_1, \ldots, g_n$ with corresponding empirical probability measure $\text{Pr}^{(n)}_G$ is given by
\begin{align*}
	D_G^{(n)}:G \times \varkappa_G \to \mathbb{R}_{\ge 0}, (g, \mathbb{K}) \mapsto D_G(g, \mathbb{K}, \text{Pr}^{(n)}_G).
\end{align*}
\end{definition}

Serving as an example, we consider the generalised Tukey depth, based on \cite{schollmeyer17b, schollmeyer17a} and introduced in \cite{blocher22}. 
The Tukey depth on $\mathbb{R}^d$ of a point $g\in \mathbb{R}^d$, c.~f.~\cite{tukey75, donoho92}, is the smallest probability of a halfspace containing $g$. To build the bridge to FCA, we consider the formal context $\Kbb$ with the object set $G = \mathbb{R}^d$, attribute set $M = \{H \subseteq \mathbb{R}^d \mid H \text{ topologically closed halfspace}\}$ and binary relation $I$ with $(g,H) \in I$ if and only if $g \in H$, see Example~\ref{exampl: spatial_1}. Based on this the Tukey depth for a point $g \in \mathbb{R}^d$ and probability measure $\Pr_{\mathbb{R}^d}$ on $\mathbb{R}^d$ can be written as
\begin{align}\label{tukey_identity}
	\inf_{H \in \Hcal(g)}  \text{Pr}_{\mathbb{R}^d}(H) = \inf_{H \in \Psi(g)}  \text{Pr}_{\mathbb{R}^d}(H) = 1 - \sup_{H \in \Psi(g)}  \text{Pr}_{\mathbb{R}^d}(G\setminus H) = 1 - \sup_{H \in M \setminus \Psi(g)}  \text{Pr}_{\mathbb{R}^d}(H)= 1 - \sup_{m \in M \setminus \Psi(g)}  \text{Pr}_{\mathbb{R}^d}(\Phi(m))
\end{align} 
where $\Hcal(g)$ is the set of all topologically closed halfspaces containing $g$ and $\Phi, \Psi$ correspond to the derivation operators given by the formal context $\mathbb{K}$ defined in Example~\ref{exampl: spatial_1}. The first and the last equality are translations into FCA language, where importantly the final right hand side does not involve the notion of a halfspace, which allows the generalisation of the Tukey depth to any arbitrary FCA setting.
The third equality holds because $G\setminus H$ is a topologically open halfspace and the supremum does not change when considering topologically open halfspaces instead. 
The right hand side of Equation (\ref{tukey_identity}) will now be the basis of the generalisation of Tukey depth to arbitrary formal contexts. Before we proceed, we shortly indicate, why we do not use directly the left hand side of Equation (\ref{tukey_identity}): Generally, the topologically closed convex sets in $\mathbb{R}^d$ correspond to the extents within our approach to use FCA to define data depth functions. On the other hand, the topologically closed halfspaces of $\mathbb{R}^d$ have no general natural equivalent in FCA. 
On the left hand side, if we replace the family of topologically closed halfspaces with the family of all extents, we obtain a depth value of zero for all $g \in \mathbb{R}^d$ when, e.g., the probability measure is continuous w.r.t.~the Lebesgue measure. This is of course unsatisfying. Moreover, note that unlike halfspaces, the complement of extents in FCA are generally not extents. This further separates the left and the right hand side of Equation~(\ref{tukey_identity}). Therefore, we take the right hand side of Equation (\ref{tukey_identity}). As will be shown later in Theorem~\ref{theo:tukey_extent} of Section~\ref{sec:tukeys_depth}, taking the supremum over all halfspaces or taking it over all topologically closed convex sets in $\mathbb{R}^d$ gives the same result, which also translates to the generalised Tukey depth where the supremum over all extents and the supremum over all extents generated by one single attribute coincide. 

Additionally, the supremum of the right hand side of Equation~(\ref{tukey_identity}) has a natural interpretation as a measure of outlyingness that can be expressed in the language of FCA, see \cite[Section~2]{schollmeyer17b} and \cite[Section~5]{schollmeyer17a}. In these articles, the author constructs special representative extents: 
For given $\alpha \in[0,1]$ we say that an extent is $\alpha$-extensive if it contains at least a proportion of  $\alpha \cdot 100\%$ of data points or probability mass. Of course, for given $\alpha$ there are many such extents but one can take for one $\alpha$ the intersection of all these extents. This intersection is again an extent which is then in a certain sense a representative extent w.r.t.~a level $\alpha$. The outlyingness of a point $g$ is then given by the (empirical) probability mass of the most specific depth contour (i.e., the most specific representative extent that corresponds to the smallest possible $\alpha$) that still contains $g$. A slightly different, but order-theoretically equivalent, definition is to take the least specific depth contour not containing $g$. This is exactly what is expressed by $\sup_{m \in M \setminus \Psi(g)}  \text{Pr}_{\mathbb{R}^d}(\Phi(m))$ in Equation~(\ref{tukey_identity}). With this motivation (and with the insight of Theorem~\ref{theo:tukey_extent}) we get as generalised Tukey depth:

\begin{definition}\label{def:gen_tukey}
Let $G$ be a set, $\varkappa_G$ a subset of formal contexts and $\Prbb_G$ a subset of probability measures on $G$. Assume that $\varkappa_G$ and $\Prbb_G$ are defined as in Definition~\ref{def:allg_depth_vorschrift}. The generalised Tukey depth is given by
\begin{align*}
	T_G: \begin{array}{l}
	G \times \varkappa_G \times \Prbb_G \to [0,1],\\
	(g, \mathbb{K}, \Pr_G) \mapsto 1 - \sup_{m \in M \setminus \Psi(g)}\Pr_G(\Phi(m)) 
	\end{array}
\end{align*}
where $\mathbb{K}$ defines the operator $\Phi$ and $\Psi$. We set $\sup_{t \in \emptyset} f(t) \defeq 0$ for every function $f$. 
\end{definition} 
Note that here, we take the supremum only over all extents generated by one single attribute, as this gives the same result as taking the supremum over all extents (c.f., Theorem~\ref{theo:tukey_extent}) and is at the same time easier to compute. Additionally, from this definition it becomes clear that the (generalised) Tukey depth does not dependent on the dependence structure between the attributes, but only on the marginal distribution of the attributes. Observe that the term marginal is meant here as the probability of lying in a certain halfspace and should not be confused with the marginal distribution of points in $\mathbb{R}^d$ in the sense of the distribution of one coordinate. In the sequel, we will always refer to the marginal distribution as the distribution of the attributes.

Furthermore, note that $\varkappa_G$ is not restricted to any subset and, in particular, $\Pr_G$ is only restricted by $\varkappa_G$, see Definition~\ref{def:allg_depth_vorschrift}. The second part of the mapping rule in Definition~\ref{def:gen_tukey}, $\sup_{m \in M \setminus \Psi(g)}\Pr_G(\Phi(m))$, corresponds to the supremum of the probabilities of the events which consists of all objects having an attribute the object of interest does not have. Thus, if $g \in G$ is an object that has all the attributes which occur, then $g$ has a maximal depth of value one. Note, however, that there are usually no objects that have all attributes. For example, consider the spatial context described in Example~\ref{exampl: spatial_1}, where the generalised Tukey depth corresponds to the well-known Tukey depth on $\mathbb{R}^d$. Thus, in this situation, the (generalised) Tukey depth is bounded by $0.5$ and, in particular, strictly below $1$. Another example is the interordinal scaling described in Section~\ref{sec:motivation}.

Now, we define the empirical version. Analogously to Definition~\ref{def:allg_empir_depth_vorschrift}, let $g_1, \ldots, g_n$ be a sample of $G$ based on $\Pr_G \in \Prbb_G$. Then the empirical generalised Tukey depth function corresponds to $T_G$ by inserting the corresponding empirical probability measure $\Pr_G^{(n)}$. This gives us:

\begin{definition}\label{def:empirical_tukey}
Let $G$, $\varkappa_G$ and $\Prbb_G$ be defined as in Definition~\ref{def:allg_empir_depth_vorschrift}. Let $(g_1,\ldots,g_n)$ be a sample from $G$ according to $\text{Pr}_G \in \Prbb_G$ with associated empirical measure $\text{Pr}_G^{(n)}$. Then the empirical generalised Tukey depth is given by
\begin{align*}
	T^{(n)}_G:
	\begin{array}{l}
	G \times \varkappa_G \to \mathbb{R}_{\ge 0},\\
	(g, \mathbb{K}) \mapsto 1 - \sup_{m \in M \setminus \Psi(g)} \sum_{\tilde{g} \in \Phi(m)}\text{Pr}_G^{(n)}(\tilde{g}) {\;= 1 - \frac{1}{n} \sup_{m \in M\setminus \Psi(g)} \# \Phi(m)}.
	\end{array}
\end{align*}
\end{definition}
As derived above, the generalised Tukey depth applied to the spatial context defined in Example~\ref{exampl: spatial_1} gives us the well-known Tukey depth on $\mathbb{R}^d$, see~\cite{tukey75}. In particular, on $\mathbb{R}^1$ we get the ranking given by the classical quantiles, with the median being the central one and the depth values decreasing outwards.
\begin{example}\label{examp:tukey_rechenbsp}
	Recall Example~\ref{motivation_titanic} and~\ref{exampl: nominal_ordinal_2}. Let us assume that the probability measure $\Pr_G$ in $G = \{g_1, g_2, g_3, g_4, g_5\}$ is uniform and the formal context $\Kbb$ is given by Table~\ref{tab:titanic_motivation}. Then the generalised Tukey depths are
	\begin{align*}
		&T_G(g_2, \Kbb, \text{Pr}_G) = 1 - \sup_{m \in M \setminus \Psi(g_2)} \text{Pr}_G(\Phi(m)) = 1- \frac{3}{5} = \frac{2}{5} = T_G(g_1, \Kbb, \text{Pr}_G) = T_G(g_5, \Kbb, \text{Pr}_G),\\
	 	&T_G(g_4, \Kbb, \text{Pr}_G) = 1 - \sup_{m \in M \setminus \Psi(g_4)} \text{Pr}_G(\Phi(m))  = 1 - \frac{4}{5} = \frac{1}{5} = T_G(g_3, \Kbb, \text{Pr}_G).
	\end{align*}
\end{example}

\begin{example}\label{examp:posets_tukey_depth}
	In Example~\ref{motivation_posets} we briefly described the formal context in which partial orders are represented. See also~\citep{blocher22, blocher23}. 
	Recall that the attributes represent the existence or non-existence of a dominance structure between two different algorithms $i$ and $j$. Thus, each pair of algorithms $i$ and $j$ gives rise to four binary attributes: First, whether $i$ dominates $j$ or vice versa, and second, whether $i$ does not dominate $j$ and vice versa.
	As can be seen from the definition, the generalised Tukey depth only considers the marginal distribution of the attributes.  
	Thus, based on the formal context describing partial orders, the generalised Tukey depth relies only on the probability of how often dominance between two classifiers occurs and how often it does not. More specifically, if we consider a
	sample, then the marginal distribution is given by the proportion of observed dominance between two classifiers, or the proportion of observed non-dominances. In the following, we omit the case where all partial orders have exactly the same probability, see~\cite{blocher22} for details on this. This assumption implies that the partial order(s) with the minimum depth value are those that do not fulfill a (non-)dominance structure that has the highest proportion. There is no such simple answer for the maximum depth partial order. This is due to the transitivity assumption on the dominance structure. A partial order with highest depth value is one that summarises the most dominance and non-dominance structures that are most often observed.
\end{example}

\begin{example}\label{examp:hierarchical_nominal_tukey}
	Recall the nominal hierarchical data structure discussed in Example~\ref{hierarchical nominal}. Note that the groups are already disjoint after the first partition and are then successively partitioned into smaller subsets. In the following we assume that the probability mass is not strongly focused on a small subset of possible observations, see the proof of Theorem~\ref{theo: tukey_strongly_free} for details. Since the generalised Tukey depth includes only the maximum marginal probability of the attributes that are not true for the object of interest, we can observe that only the first partition is used for the generalised Tukey depth. Therefore, dropping all other attributes (representing the finer splits) and thus considering a smaller and less informative formal context does not change the generalised Tukey depth based on this nominal hierarchical formal context. In other words, the generalised Tukey depth ignores the further information of the data.
\end{example}
Before presenting the structural properties, we define when two depth functions are isomorph and thus represent the same center-outward order.
\begin{definition}\label{def:isomorph}
	Let $G$ be a set and $D_G$ and $\tilde{D}_G$ be two depth functions based on $\varkappa_G$ and $\Prbb_G$, $\tilde{\varkappa_{G}}$ and $\tilde{\Prbb}_{G}$ respectively, see Definition~\ref{def:allg_depth_vorschrift}. Let $\Kbb \in \varkappa_G$, $\tilde{\Kbb} \in \tilde{\varkappa}_{G}$, $\Pr_G \in \Prbb_G$ and $\tilde{\Pr}_G \in \tilde{\Prbb}_{G}$. Then $D_G(\cdot, \Kbb, \Pr_G)$ and $\tilde{D}_G(\cdot, \tilde{\Kbb}, \tilde{\Pr}_{G})$ on $G$ are \textit{isomorph} if and only if there exists a bijective and bimeasureable function 
	$i: G \to G$ such that 
	\begin{align*}
		D_G(g, \Kbb, \text{Pr}_G) \le D_G(\tilde{g}, \Kbb, \text{Pr}_G) \Longleftrightarrow \tilde{D}_G(i(g), \tilde{\Kbb}, \tilde{\Pr}_{G}) \le \tilde{D}_G(i(\tilde{g}), \tilde{\Kbb}, \tilde{\Pr}_{G}) 
	\end{align*}
	is true for all $g, \tilde{g} \in G$. We call a bijective function \textit{bimeasurable} if and only if $i$ and $i^{-1}$ are measurable w.r.t.~the corresponding $\sigma$-fields. In what follows, $\cong$ denotes the isomorphism between two depth functions. 
\end{definition}

For simplicity of notation, we write $D$, $D^{(n)}$, $T$, $T^{(n)}$, $\varkappa$, $\Prbb$, $\Pr$ and $\Pr^{(n)}$ instead of $D_G$, $D_G^{(n)}$, $T_G$, $T_G^{(n)}$, $\varkappa_G$, $\Prbb_G$, $\Pr_G$ and $\Pr_G^{(n)}$ in the following if the underlying object set $G$ is clear.

\section{Structural properties characterising the depth function using formal concept analysis}\label{sec:des_properties}
\begin{figure}
\begin{center}
	\includegraphics[trim={0 0 0 0},clip, width=0.6\textwidth]{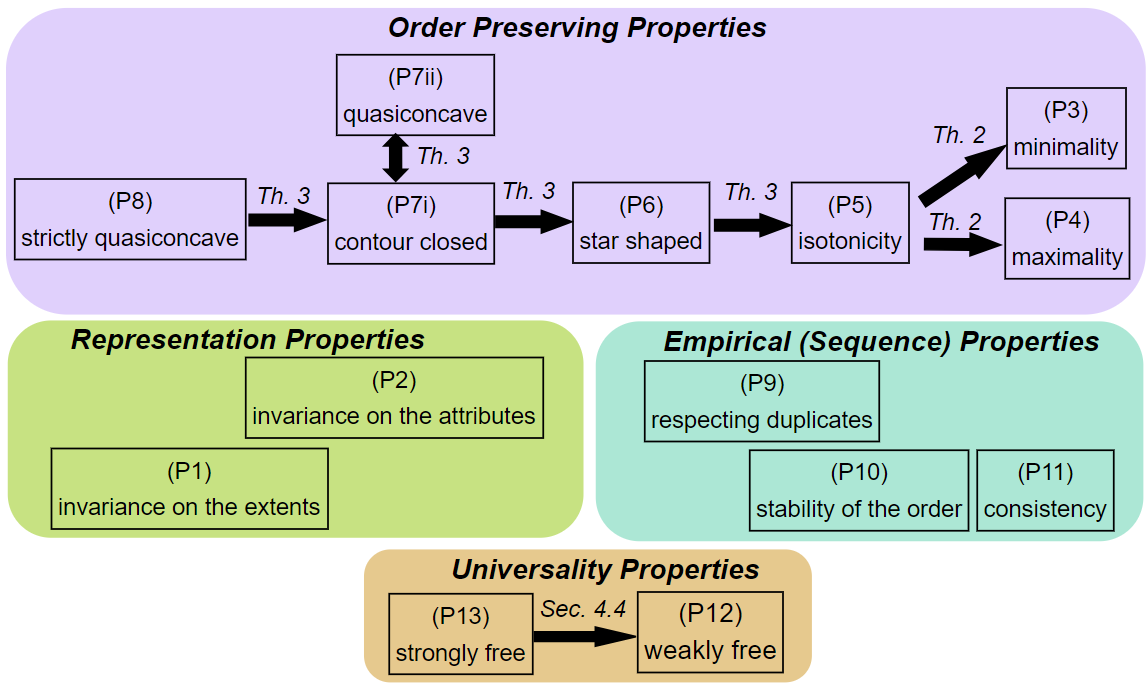}
\caption{\label{fig:uebersicht_property_structure} Overview of the structural properties together with their mathematical connections. \textit{Th} is Theorem and \textit{Sec} is Section for short. This overview contains all direct implications between the properties. Further results like, e.g., limitations of the properties, see, e.g., Theorem~\ref{prop:strength_strc_conc}, are not presented in the Figure.}
\end{center}
\end{figure}

After this general definition of depth functions based on FCA, we want to discuss some \textit{structural properties} a depth function can have. With this, we tackle the question of what centrality is and which object is - in some sense - \textit{closer} to the center than another object. In particular, we provide concepts to discuss centrality and outlyingness for different data types without necessarily presupposing one specific data structure. Thus, this section gives a starting point for a discussion on centrality, outlyingness and depth functions based on FCA. Furthermore, we define a framework upon which newly introduced depth functions can be studied and compared.

For normed vector spaces, there is already an ongoing discussion about this. The authors of \citep{liu90, zuo00, mosler13, mosler22} are concerned about these questions for data depth functions defined on $\mathbb{R}^d$. Furthermore, \cite{gijbels17, durocher22} discuss depth functions and centrality topics for functional data. All these examples have in common that they are based on a normed vector space and that the clarification of centrality and outlyingness is done by defining properties. For example, \cite{zuo00, mosler13} demand that the depth of a point $x\in\mathbb{R}^d$ should converge to zero as the norm of $x$, i.e. $||x||$, tends to infinity. This reflects the intuition about outlyingness in unbounded spaces. In contrast to outlyingness, the definition of a center point does not immediately follow. Basically, every point in $\mathbb{R}^d$ can be the center. This follows from the fact that after translation, rotation, etc. the structure of a normed vector space does not change. The center generally seems to be only naturally specified in special cases. For example \cite{zuo00} assumes that for a probability measure that is symmetric around a point $c$, see \cite{zuo00b}, this point $c$ should then be the center. In the case of a multivariate normal distribution, the center point is the mean vector. Since a depth function should represent a center-outward order, this point $c$ needs to have a maximal depth value.
Together with the center, one needs to discuss what it means that one point is further away from the center than another.
In $\mathbb{R}^d$, this is achieved by the use of line segments, see \cite{zuo00, mosler13}. A point $p_1$ which lies on the line segment between the center and a further point $p_2$ is said to be closer to the center than $p_2$. In other words, $p_2$ is more outlying than $p_1$. At the first glance, here we also use the normed vector space structure.

It follows from the above that the definition of the terms centrality and outlyingness seem to highly rely on the underlying distribution. In particular, the part defining the center point which is based on a symmetric distribution, see \cite{zuo00b}. 
Slightly different to \cite{zuo00,zuo00b} is the approach of \cite{mosler13}. 
The properties given by \cite{mosler13} put emphasis on the structure of the underlying spaces and not on the probability measure. Thus, in contrast to e.g., \cite{zuo00,zuo00b}, the definition of centrality of a point is more detached from the notion of symmetry/centrality of the underlying probability measure. Note that \cite{mosler13} gives an even stronger definition of depth functions. In the following, we discuss the term centrality and outlyingness from the perspective of the approach in \cite{mosler13}.

In this article, we consider a space $G$ where the structure of the data points is given by a formal context. In the style of \cite{zuo00} and \citep{mosler13} for depth functions in $\mathbb{R}^d$, we define \textit{structural properties} of a depth function based on FCA to clarify the notion of centrality and outlyingness. These properties express characteristics of the data set when represented via a formal context and how these characteristics are included in the depth function. Some of the structural properties build upon the existing ones in $\mathbb{R}^d$ that we transfer to our new situation. To achieve this, we use that the set of all topologically closed convex sets is a natural closure system on $\mathbb{R}^d$. Since the extent set also defines a closure system, we obtain a direct translation to FCA by representing the properties by use of the convex closure system on $\mathbb{R}^d$. For example, we transfer the idea of a line segment to our data representation via a formal context, see Property (P6). Moreover, the \say{quasiconcavity} property, see \cite[p. 19]{mosler13} has a natural translation to FCA by use of the closure system, see Property (P7i and ii). Other structural properties introduced discuss the possibility that the data itself may have a center-outward order. This inherited order is then represented by the formal context. For example, consider the extreme case that an object lies in every set of the closure system. This means that it has every attribute that is true for any element of the data set. Thus, this object is implied by every other object and is therefore more specific than any other object. So it should have maximum depth for any possible probability measure, see Property~(P4). Of course, the other extreme case, where an object lies only in the extent corresponding to the entire set and therefore shares no attributes with the other data elements, should have minimal depth by default, see Property~(P3). Note that these two extremes do not necessarily occur. For example, there is no such element in $\mathbb{R}^d$ together with the convex sets. All in all, one can say that the introduced structural properties are defined and discussed within two perspectives: The transfer of already existing properties in $\mathbb{R}^d$ and the new development of properties based on the theory of FCA, reflecting the inherited center-outward structure of the data.

In total 13 structural properties are presented which can be covered under four different categories: Representation properties, order preserving properties, empirical (sequence) properties and universality properties. An overview of the structural properties together with their mathematical connections can be seen in Figure~\ref{fig:uebersicht_property_structure}. In what follows, we fix the set of objects $G$ and consider the depth function $D: G \times \varkappa \times \Prbb\to \mathbb{R}_{\ge 0}$. Furthermore, we refer to Definition~\ref{def:allg_depth_vorschrift} when we say that $D$ is a depth function. We use the term empirical depth function when considering a depth function in the style of Definition~\ref{def:allg_empir_depth_vorschrift}. 
Afterwards, in Section~\ref{sec:tukeys_depth} we check these properties based on the introduced generalised Tukey depth, see Section~\ref{sec:definition}.

\subsection{Representation properties}\label{sec:representation_properties}
Depth functions satisfying the following properties are structure preserving on $G$. This includes two parts: First, assume that we represent the data set by two formal contexts and that we consider two probability measures. Let us assume that between the corresponding extent sets and therefore the structure of the objects exists a bijective function. Furthermore, if additionally the probability values are preserved by this bijective function, then the depth function should be preserved as well. In other words, the order given by the depth function should not rely on the used scaling method unless it does substantially change the structure of the extent sets. Also, the influence of the probability measure is only based on the extent sets. In Section~\ref{sec:motivation} we motivated this property by illustrating that the observation values are represented within the arrangement of the different groups of the passengers/observations. Thus, if neither the groups nor the probabilities of these groups change, then the depth value should also remain the same. This can be seen as an adaptation of the \say{affine invariance} property of depth functions in $\mathbb{R}^d$, see \cite[p. 463]{zuo00}. There, the depth function equals the depth of the shifted version if the probability measure is shifted accordingly.

The second part of the representation properties considers the attribute values. Here, we say that the depth function preserves the structure if and only if two objects with the same attributes have the same depth value. From the perspective of FCA, two objects with the same attributes are duplicates. Therefore they should be assigned to the same depth value. In the formal context of the Titanic example, see Example~\ref{motivation_titanic}, this means that two passengers of exactly the same age, sex and class of travel must have the same depth values.

\begin{enumerate}
\item[(P1)] \textit{Invariance on the extents:} Let $\Kbb, \tilde{\Kbb} \in \varkappa$ be two formal contexts on $G$ and let $\Pr, \tilde{\Pr} \in \Prbb$ be two probability measures on $G$. If there exists a bijective and bimeasureable function $i: G \to G$ such that the extents are preserved (i.e.~$E$ extent w.r.t.~$\Kbb$ $\Leftrightarrow i(E)$ extent w.r.t.~$\tilde{\Kbb}$) and the probability is also preserved (i.e. $\Pr(E) = \tilde{\Pr}(i(E))$), then
	\begin{align*}
		D(\cdot, \Kbb, \Pr) \cong D(\cdot, \tilde{\Kbb}, \tilde{\Pr})
	\end{align*}
	is true.
\end{enumerate}

\begin{enumerate}
	\item[(P2)] \textit{Invariance on the attributes:} For every $\Kbb \in \varkappa$, $\Pr \in \Prbb$ and $g_1, g_2 \in G$ with $\Psi(\{g_1\}) = \Psi(\{g_2\})$, 
	\begin{align*}
		D(g_1, \Kbb, \Pr) = D(g_2, \Kbb, \Pr)
	\end{align*}
	holds.
\end{enumerate}

\begin{remark}
	Property~(P1) is the only property that directly discusses the connection between two different formal contexts. This is due to the fact that even just adding more attributes can strongly change the implications and the whole data structure representation, even though this procedure \say{only} adds more extent sets. Therefore, in general, it is not clear how the depth function should behave when more attributes are added. However, for example for interordinal scaling, which is \say{finer} than ordinal scaling, see Example~\ref{motivation_titanic}, it can be shown that the following Property~(P7i and~ii) together with Property~(P1) in many cases already cover the desired changes in the order of the data elements. More precisely, for ordinal scaling, Property~(P7i and~ii) ensures that the smallest value has maximum depth, while for interordinal scaling it gives more flexibility and, depending on the underlying probability measure, any element can have the highest depth value. In particular, in many cases (e.g. uniform distribution) Properties~(P1) and (P7i and~ii) imply that the classical median has the highest depth. For the generalised Tukey depth, we always obtain the classical median as the highest depth for interordinal scaling.
\end{remark}

\subsection{Order preserving properties}\label{sec:order_preserving_properties}

While the representation properties ensure that similar structures on $G$ lead to the same depth function, the next \textit{order preserving properties} consider the obtained order by the depth function. These properties increase in their strength of restriction.

These properties are defined along the lines of \say{monotonicity relative to the deepest point} and \say{maximality at the center} properties, see \cite[p. 463]{zuo00}, and the \say{quasiconcavity} properties, see \cite[p. 19]{mosler13}, defined for $\mathbb{R}^d$. In contrast to $\mathbb{R}^d$ we neither have a norm nor a concept of translation and symmetry, but we can make use of the closure system and closure operator given by the extent set. 
When considering the properties defined for $\mathbb{R}^d$ in the context of the convex sets which define a closure system, we can build a bridge to FCA. In particular, the \say{quasiconcavity} property, see \cite[p. 19]{mosler13}, has therefore a natural adaptation. Nevertheless, the convex sets are a special case of closure sets: For example, the affine invariance is reflected in the set since by shifting the set of convex sets the family of sets does not change. This does not hold in general for closure systems. For example, there exist closure systems where one point occurs more often in the sets of the closure system than another point.
We use the concept of a formal context where two natural extreme opposite characteristics of the objects can occur.  The first one: If an object has every attribute then it lies in every extent set. Conversely, if an object has no attribute at all, then the only extent containing this object is the entire set $G$.
Property (P3) and (P4) now ensure that these two opposite characteristics are also reflected in the depth function. 
\begin{enumerate}
	\item[(P3)] \textit{Minimality:} 
	Let $\Kbb \in \varkappa, \Pr \in \Prbb$. Further, let $g_{non}\in G$ such that for every extent $E\subsetneq G$ of $\Kbb$ 
	we have that $g_{non} \not\in E$, then
	\begin{align*}
		 D(g_{non},\Kbb, \Pr) = \min_{g\in G} D(g, \Kbb, \Pr)
	\end{align*}	
	is true.
	
	\item[(P4)] \textit{Maximality:} Let $\Kbb \in \varkappa, \Pr \in \Prbb$. Assume there exists $g_{all} \in G$ such that for every extent $E$ of $\Kbb$ we have that $g_{all} \in E$. Then $$D(g_{all}, \Kbb, \Pr) = \max_{g \in G} D(g, \Kbb, \Pr)$$ holds.
\end{enumerate}

Note that a depth function which fulfils these two properties does not rely on the probability measure to set the values of the two extreme cases. On the other hand, there exist formal contexts such that the objects $g_{non}$ or $g_{all}$ do not exist at all. Recall Example~\ref{exampl: spatial_1} where we considered the spatial data. 

Nevertheless, since we have now predefined the maximal depth value in specific cases, this has to be in line with adapting properties like ``monotone on rays'' or ``quasiconcavity'', see \cite[p. 19]{mosler13}. In what follows we start with less restricting properties and increase their restriction strength slowly. We show that they imply Properties (P4) and (P3).

Property (P5) is inspired by the fact that in FCA, an object $g_2$ which lies in the closure of an other object $g_1$ implies that this object is more specific than $g_2$. In other words, $g_2$ has all attributes and possibly even more attributes than $g_1$. This is analogous to the assumption that $\gamma(\{g_1\}) \supseteq \gamma(\{g_2\})$ is true. Thus, Property (P5) says that an object $g_2$ must have a depth value as least as high as $g_1$. 
\begin{enumerate}	
	\item[(P5)] \textit{Isotonicity}: For every $\Pr \in \Prbb$ and formal context $\Kbb \in \varkappa$ with $g_1, g_2 \in G$ such that $\gamma_{\mathbb{K}}(\{g_1\}) \supseteq \gamma_{\mathbb{K}}(\{g_2\})$,	\begin{align*}
		D(g_1, \Kbb, \Pr) \le D(g_2, \Kbb, \Pr) 
	\end{align*}
	is fulfilled.
\end{enumerate}
There exists a natural strengthening of the isotonicity property (P5) which leads to the adaptation of the ``monotonicity relative to deepest point'', see \cite[p. 463]{zuo00}, property in $\mathbb{R}^d$. For start, let us assume that the depth function is bounded from above. Furthermore, we assume that the depth function has its maximum at center $c \in G$. The ``monotone on rays'' in $\mathbb{R}^d$ definition in \cite[p. 19]{mosler13} states that the depth of a point that moves further away from the center $c$ on a fixed ray should decrease. Since we do neither have a norm nor a vector space, we cannot generally define what \textit{further away} as well as \textit{ray} means. Thus, we translate the definition to a setting using the convex closure operator instead. In this case, if a point/object $\tilde{g}$ lies on the line segment given by a further point $g$ and the center $c$, then the depth of $\tilde{g}$ must be at least as high as the depth of $g$. In other words, when a point $\tilde{g}$ lies in the convex closure of another point $g$ and the center, then we have a lower bound for the depth of $\tilde{g}$. Thus, the points/objects with depth values larger or equal to a fixed value form a starshaped set. This gives the name of the property. With the definition based on the convex closure system, we can easily transfer this property to FCA and obtain Property~(P6). 
\begin{enumerate}	
	\item[(P6)] \textit{Starshapedness}: Let $\Kbb\in \varkappa$ and $\Pr \in \Prbb$. If there exists at least one center point $c \in G$ such that for all $g \in G$ and all $\tilde{g} \in \gamma_{\mathbb{K}}(\{c,g\})$ we have $$D(\tilde{g}, \mathbb{K}, \Pr) \ge D(g, \mathbb{K}, \Pr),$$ then we call $D$ starshaped.
\end{enumerate}

In the above explanation, we assumed that the center point $c$ has a maximal depth value. This assumption is not included in the definition of Property (P6) as the boundedness of the depth function, as well as the maximality at the center, follow directly. Let us fix $c\in G$ to be one center point. Since for every closure operator the isotonicity assumption (see Section~\ref{sec:preliminaries}, not to be confused with Property (P5)) holds, we get that for every object $g$ we have $c \in \gamma(\{c,g\})$. Thus, by Property (P6) for every $g \in G$ we obtain
\begin{align*}
	D(c, \mathbb{K}, \Pr) \ge D(g, \mathbb{K}, \Pr).
\end{align*} 
With this, $c$ must have the highest depth value. Since the depth function maps to $\mathbb{R}$, this gives us also the upper bound of the depth function. Note that a center $c$ is not always naturally given. Besides the spatial data, see Example~\ref{exampl: spatial_1}, where one can use the probability measure to get a center, in the framework of FCA one can say that $g_{all}$ in Property~(P4) could be such a center. In contrast, in the example of the Titanic data, see Table~\ref{tab:titanic_motivation}, there is no passenger who is naturally a center. Each passenger has at least one component where their observation is not in the center (median for age and modus for sex and class of travel).

 Furthermore, the starshaped property (P6) together with the invariance on the extents property (P1) has some implications on the point of maximal depth when the underlying probability measure has some symmetry property. 
As already indicated at the beginning of Section~\ref{sec:des_properties}, because of the lack of a translation operation, etc., it is difficult to define symmetry in our setting.
However, one can still define some notion of point symmetry, see the next theorem.

\begin{theorem}\label{theo: symmetry based on self-inverse}
    Let $\Kbb \in \varkappa$ and $\Pr \in \Prbb$. We assume that $\Pr$ is point symmetric around $s \in G$ in the following sense:
    There exists a bimeasurable involutory function $i:G\longrightarrow G$ (i.e. $i(i(g)) = g$ for all $g \in G$) such that
    \begin{enumerate}
    		\item[1.] $\Pr(i(E))=\Pr(E)$ for all extents $E\subseteq G$ and
    		\item[2.] for all $g \in G$ we have $s \in \gamma(\{g, i(g)\})$.
    \end{enumerate}
Then for every depth function which fulfils Properties (P1) and (P6), the center of symmetry $s$ of $\Pr$ has maximal depth. 
\end{theorem}

\begin{proof}
 Let $c$ be one center point. First note that by Assumption 2 we get $s \in \gamma (\{c, i(c)\})$. Because $D$ is starshaped, we have $D(s,\Kbb,\Pr) \geq 
    D(i(c),\Kbb,\Pr)$. Now, we use that there exists a bimeasurable involutory function $i$ and that $D$ is invariant on the extents (P1) to get $D(i(c),\Kbb,\Pr)= D(c, \Kbb, \Pr)$. With this, we conclude that $D(s,\Kbb,\Pr)\geq D(c,\Kbb,Pr)$. This proves that $s$ has maximal depth.
\end{proof}

Note that object $s$ of Theorem~\ref{theo: symmetry based on self-inverse} is not necessarily a center point in the style of Property (P6).

These properties are written down in their order of strength. More precisely, Property (P6) implies Property (P5) and Property (P5) implies Properties (P4) and (P3), see Theorem~\ref{lem:P5_implies_P3_4}. Thus, if Property (P6) is satisfied and if there exists an object $g$ which has every attribute, then $c$ must be one of the center point $c$ discussed in Property (P6).
\begin{theorem}\label{lem:P5_implies_P3_4}
Let $\Kbb \in \varkappa$ and $\Pr \in \Prbb$. Let $D$ be a depth function then the following implications hold for $D$:
\begin{enumerate}
	\item Let $c \in G$ be a center object. If $D$ satisfies Property (P6), then Property (P5) is true for $D$.
	\item If $D$ satisfies Property (P5), then Properties (P4) and (P3) are true for $D$.
\end{enumerate}
\end{theorem}
\begin{proof}
Assume that $\Kbb \in \varkappa$ and $\Pr \in \Prbb$. Let $D$ be a depth function.

We begin by proving that Property (P6) implies Property (P5). Therefore assume that Property (P6) is true for $D$ and let $c \in G$ be one center point. Further let $g_1, g_2 \in G$ such that $\gamma_{\mathbb{K}}(\{g_1\}) \supseteq \gamma_{\mathbb{K}}(\{g_2\})$. Then we get that $g_2 \in \gamma(\{c, g_1\})$ and therefore we have $D(g_1, \Kbb, \Pr) \le D(g_2, \Kbb, \Pr)$.

The next step is to show that Property (P5) implies  Property (P4). Assume that Property (P5) is satisfied and that $g_{all}$ lies in every extent set. Then, we get for every $g \in G$ that $\gamma_{\mathbb{K}}(\{g_{all}\}) \supseteq \gamma_{\mathbb{K}}(\{g\})$ holds. Thus, $D(g_{all}, \Kbb, \Pr) \ge \max_{\{g \in G\}} D(g, \Kbb, \Pr)$ is true and Property (P4) follows.

Finally, we show that Property (P3) follows from Property (P5). Let $g_{non}\in G$ be an object which lies only in the entire set and in no other extent set. Furthermore, we suppose that Property (P5) is true for $D$. Since $\gamma(\{g_{non}\}) = G$ for every $g \in G$ we have $\gamma(\{g\}) \subseteq \gamma(\{g_{non}\})$. Due to Property (P5) we follow that $D(g_{non}, \Kbb, \Pr) \le D(g, \Kbb, \Pr)$ for every $g \in G$. This gives us Property (P3).
\end{proof}

\begin{remark}\label{rem:center point}
Let us point out some consequences of Property (P6). In contrast to \cite[p. 463]{zuo00} we do not assume that the center $c$ is unique. For example, the depth function is allowed to have a plateau at the highest point. In particular, we allow the depth function to be constant.  Especially $g_{non}$ can only be the center point if the function is constant, due to Properties (P3), (P4) and Theorem~\ref{lem:P5_implies_P3_4}. Moreover, we get that when the depth function has at least two different values, then $g_{non}$ must have the minimal value and $g_{all}$ the maximal value. This allows us to specify a center point and an outlying point without relying on the other observed points/objects by defining the scaling method. Note that this situation does not occur often. For example, all scaling methods described in Section~\ref{sec:preliminaries} do not have an object $g_{non}$ and~$g_{all}$.

Nevertheless, this stresses the importance of a meaningful and carefully chosen scaling method. Observe that the difference between having an attribute and not having an attribute is not symmetric and cannot be switched without eventually fundamentally changing the characteristic of the corresponding closure system.
\end{remark}

The next order preserving properties are in the style of the \say{quasiconcavity} property, see, e.g., \cite[p. 19]{mosler13}. In the context of $\Rbb^d$ quasiconcavity states that the set of all points with a larger (or equal) depth value than $\alpha \ge 0$ needs to be a convex set. These sets are called contour sets. In this case, there is a direct transfer to FCA by the extent set as a closure system. Therefore, we first have to define the contour sets within the theory of FCA. Let $\Kbb\in \varkappa$ be a formal context and $\Pr \in \Prbb$ be a probability measure. For $\alpha \in  \text{im}(D(\cdot, \Kbb, \Pr))$ the contour set $Cont_{\alpha}$ is defined as follows
\begin{align*}
	Cont_{\alpha} \defeq \{g \in G \mid D(g, \Kbb, \Pr) \ge \alpha\}.
\end{align*}
Now, we say the depth function is countourclosed if every contour set is an extent set. This is stated in Property (P7i). Instead of considering the contour sets, an analogous statement for Property (P7i) is given in Property (P7ii). Here, we assure that the depth of an observation that lies in the closure of an input set is larger or equal to the infimum of the input. 
Property (P7ii) has a natural strengthening by assuming strict inequalities instead, see Property (P8).

\begin{enumerate}
	\item[(P7i)] \textit{Countourclosed}: For every formal context $\mathbb{K} \in \varkappa$, probability measure $\Pr\in\Prbb$ and every $\alpha \in  \text{im}(D(\cdot, \Kbb, \Pr))$ the contour set $Cont_{\alpha}$ is an extent of the formal context $\Kbb$.
	
	\item[(P7ii)] \textit{Quasiconcave}: Let $\Kbb\in \varkappa$ and $\Pr \in \Prbb$. If for all $A \subseteq G$ and all $g \in \gamma_{\mathbb{K}}(A) \setminus A$ we have $$D(g, \mathbb{K}, \Pr) \ge \inf_{\tilde{g} \in A} D(\tilde{g},\mathbb{K}, \Pr),$$ we call $D$ quasiconcave.
	
	\item[(P8)] \textit{Strictly quasiconcave}: Let $\Kbb\in \varkappa$ and $\Pr \in \Prbb$. If for all $A \subseteq G$ and all $g \in \gamma_{\mathbb{K}}(A) \setminus A$ we have $$D(g, \mathbb{K}, \Pr) > \inf_{\tilde{g} \in A} D(\tilde{g},\mathbb{K}, \Pr),$$ $D$ is strictly quasiconcave.

\end{enumerate}

Recall Example~\ref{motivation_titanic}. We showed that a group containing $g_3$ and $g_4$ must also contain $g_2$. Property~(P7i and ii) now formalises the statement in Example~\ref{motivation_titanic} that the depth of $g_2$ is at least as high as the minimum depth of $g_3$ and $g_4$. In Example~\ref{examp:tukey_rechenbsp} we saw that this is true for the generalised Tukey depth. In the next section we will prove that Property~(P7i and ii) holds in general for the generalised Tukey depth, but that this is not true for Property~(P8).

Again these properties have mathematical connections to each other. First of all, Property (P7i) and (P7ii) are equivalent. Secondly, Property (P8) is indeed stronger and implies Property (P7i and ii). Thirdly, if $c$ is a center point with maximal depth, then with this $c$ we can show that (P7i and ii) imply the starshaped property (P6). Moreover, Property (P7ii) implies Property (P5).

\begin{theorem}\label{theo:P7i_ii_equi_P8_impl_P7}
	 Let $\Kbb\in\varkappa$ be a formal context and $\Pr\in\Prbb$ a probability measure. Let $D(\cdot, \Kbb,\Pr)$ be a depth function. 
	 \begin{enumerate}
	 	\item Statement (P7i) and (P7ii) are equivalent.
	 	\item Property (P8) implies (P7ii).
	 	\item If their exists $c \in G$ with maximal depth value, then Property (P7ii) implies (P6) with $c$ being a center object.
	 	\item Property (P7ii) implies (P5).
	 \end{enumerate}
\end{theorem}

\begin{proof}
Let $\Kbb\in\varkappa$ and $\Pr\in\Prbb$. Let $D(\cdot, \Kbb,\Pr)$ a depth function. Note that the claim that Property (P8) implies (P7ii) is given immediately. Thus, we only show Part 1., 3. and 4. of Theorem~\ref{theo:P7i_ii_equi_P8_impl_P7}.
\begin{enumerate}
	\item[1.]First assume that Property (P7ii) is true and let $\alpha\in \text{im}(D(\cdot, \Kbb, \Pr))$ be arbitrary. We prove that Property (P7i) follows. Assume, by contradiction that $Cont_\alpha$ is not an extent set. Since $\gamma$ is a closure operator and therefore idempotent, there exists $g \in \gamma(Cont_{\alpha})\setminus Cont_{\alpha}$. Since Property (P7ii) is true, $D(g, \Kbb, \Pr) \ge \inf_{\tilde{g} \in Cont_{\alpha}}D(\tilde{g}, \Kbb, \Pr) \ge \alpha$ holds. This contradicts $g \notin Cont_{\alpha}$ and we get that Property (P7i) is fulfilled.
	
	For the reverse let (P7i) be true and let $A \subseteq G$ be arbitrary. We set $\alpha = \inf_{\tilde{g} \in A} D(\tilde{g},\Kbb, \Pr)$ and we know that $A \subseteq Cont_{\alpha}$. By (P7i) we know that $Cont_{\alpha}$ is an extent set. Since $\gamma(A)$ is the smallest extent set containing $A$ and the set of extents is a closure system, we follow that $\gamma(A) \subseteq Cont_{\alpha}$. Thus the depth of every object $g \in \gamma(A)$ must be larger or equal to $\alpha$ which implies (P7ii).
	\item[3.]Now, we assume that Property (P7ii) is true and we show that Property (P6) holds as well. Therefore, we assume that $c \in G$ has maximal depth value. Let $g, \tilde{g} \in G$ such that $\tilde{g} \in \gamma(\{c,g\})$ is true. Due to Property (P7ii) we get $D(\tilde{g}, \Kbb, \Pr) \ge \min\{D(c, \Kbb, \Pr), D(g, \Kbb, \Pr)\} \ge D(g, \Kbb, \Pr).$ The last inequality follows from the assumption that $c$ has maximal depth value. This gives us Property (P6) with $c$ being one center object.
	\item[4.] Finally, assume that Property (P7ii) holds. Let $g_1, g_2 \in G$ with $g_1 \neq g_2$ such that $\gamma(g_2) \supseteq \gamma(g_1)$. Then the quasiconcavity property implies that $D(g_1, \Kbb, \Pr) \ge D(g_2, \Kbb, \Pr)$. This shows Property (P5).
\end{enumerate}
\end{proof}


Property (P8), the strictly quasiconcavity, is indeed a very strong assumption. In particular, there exist formal contexts such that Property (P8) can never be fulfilled. An example can be found in \cite{blocher22} and Example~\ref{motivation_posets}. There, we discussed the special case of depth functions for partial orders and analysed Properties (P7i and ii) and (P8) in this case. We showed that the quasiconcavity cannot be fulfilled by the formal context of partial orders. Another example is given by Table~\ref{tab:limits_tukey_allP8} (left). Here, let $G = \{g_1, g_2,g_3\}$ and let $\Pr$ be an arbitrary probability measure. In contradiction, let us assume that Property~(P8) is true for a depth function $D$ on $G$. We have $\gamma(\{g_i, g_j\}) = \{g_1, g_2, g_3\}$ for all $i,j \in\{1,2,3\}$ with $i \neq j$. With this, we obtain that $\min\{D(g_1, \Kbb, \Pr), D(g_2, \Kbb, \Pr)\} < D(g_3, \Kbb, \Pr)$. W.l.o.g.~let $D(g_1, \Kbb, \Pr) < D(g_3, \Kbb, \Pr)$. Together with $\min\{D(g_2, \Kbb, \Pr), D(g_3, \Kbb, \Pr)\} < D(g_1, \Kbb, \Pr)$ due to Property~(P8), we have $D(g_2, \Kbb, \Pr) < D(g_1, \Kbb, \Pr)< D(g_3, \Kbb, \Pr)$. However, this is a contradiction to $\min\{D(g_1, \Kbb, \Pr), D(g_3, \Kbb, \Pr)\} < D(g_2, \Kbb, \Pr)$. Hence, Property~(P8) cannot hold for this context. In contrast, the quasiconcavity property (P7ii) allows equality which solves such problems. Note that in the case of the formal context given by Table~\ref{tab:limits_tukey_allP8} (left), Property (P7ii) leads to a constant depth function.

More generally, one can say that the strictly quasiconcavity assumption on $D$ can never be satisfied if one formal context $\Kbb \in \varkappa$ is an element of the following set:
\begin{align*}
  \Cscr^{\cancel{P8}} = \left\{(G,M,I) \text{ formal context }\bigg| \begin{array}{l} \text{exist } A, \tilde{A}\subseteq G \text{ such that } \#A<\infty, \#\tilde{A}<\infty \text{ and } A\cap\tilde{A}=\emptyset\\
  \text{and } A \subseteq \gamma(\tilde{A}), \: \tilde{A}\subseteq \gamma(A)
  \end{array} \right\}.
\end{align*}

Theorem~\ref{prop:strength_strc_conc} proves that for every $\Kbb \in \Cscr^{\cancel{P8}}$ there exists no strictly quasiconcave depth function. Here, we get a contradiction with the strict larger assumption. This case occurs naturally when the scaling method assigns to two different objects the same attribute values. 
\begin{theorem}\label{prop:strength_strc_conc}
	For every $\Kbb \in \Cscr^{\cancel{P8}}$ and every $\Pr \in \Prbb$ there exists no depth function $D$ such that Property (P8) is fulfilled.
\end{theorem} 
\begin{proof}
	Let $\Kbb \in \Cscr^{\cancel{P8}}$. We assume that there exists a depth function $D$ and a probability measure $\Pr$ such that Property~(P8) is satisfied. Since $A \subseteq \gamma(\tilde{A})$ and $\tilde{A}\subseteq \gamma(A)$ we get $\forall a \in A\setminus\tilde{A}\colon D(a, \Kbb, \Pr) > \inf_{\tilde{a}\in \tilde{A}}D(\tilde{a}, \Kbb, \Pr)$ and $\forall \tilde{a} \in \tilde{A}\setminus A\colon \inf_{a \in A}D(a, \Kbb, \Pr) < D(\tilde{a}, \Kbb, \Pr).$ Furthermore, we assumed that $\#A<\infty, \#\tilde{A}<\infty$ is true. Thus, the infimum is attained in $A$ and $\tilde{A}$. Let $a_m \in A$ be an argument of the minimum of $D(a, \Kbb, \Pr)$, and analogously we set $\tilde{a}_m \in \tilde{A}$. Since $A \cap \tilde{A} = \emptyset$, we obtain that $\tilde{a}_m\in \tilde{A}\setminus A$ and $a_m \not\in A\setminus\tilde{A}$. But with this we get $D(a_m, \Kbb, \Pr) > D(\tilde{a}_m, \Kbb, \Pr)$ and $D(a_m, \Kbb, \Pr) < D(\tilde{a}_m, \Kbb, \Pr).$ This cannot be true which contradicts the assumption of Property (P8) being true.
\end{proof}

We want to point out that $A$ and $\tilde{A}$ being finite is crucial since in the proof we used that the infimum is attained and not only the largest lower bound of the sets. Secondly, the intersection of $A$ and $\tilde{A}$ being empty is also necessary since else one can set the elements in the intersection to the minimal depth value. This is then following Property~(P8).

While the importance of a meaningful scaling method for Properties (P3) and (P4) is easily seen, this comes also into account for the last properties. For example Property (P7ii) and (P8) state that the depth of an object $g$ which is implied by a set $A$, i.e. $g \in \gamma(A)$, must have larger (or equal) depth than the minimal depth value of the elements in $A$. In the context of FCA, one can say that $g$ contains all characteristics/similarities having the objects in $A$ in common or even more. Therefore it is at least as specific (or more) than all the elements in $A$ together. Thus, when defining the scaling method, not only the individual attribute values should be taken into account, but also which object combination is more specific than others.

All these order preserving properties aim at representing the underlying data structure in the depth function. This includes structures that are also valid for $\mathbb{R}^d$, e.g. quasiconcavity, but also the extreme cases where a natural center or outlying point is given by the data structure itself, see maximality and minimality property. This emphasises that non-standard data can be very different from $\mathbb{R}^d$, and these properties aim to preserve the different structures. In the next section, we focus on the set of probability measures $\Prbb$ directly.

\subsection{Empirical (sequence) properties}\label{sec:empirical_properties}
In this section, we consider the reverse of the above section. For a fixed formal context $\Kbb$ we are interested in the behaviour of the depth function when the (empirical) probability measure changes. In what follows, we regard different empirical probability measures $\Pr^{(n)}$ induced by different samples $g_1, \ldots, g_n \in G$ with $n \in \mathbb{N}$. 
We assume that $\Pr^{(n)} \in \Prbb$ holds. Hence, in this section, we consider the empirical depth function, recall Definition~\ref{def:allg_empir_depth_vorschrift}. 

The first two properties discuss how two empirical depth functions differ if the two corresponding samples differ in a specific manner. The first one considers the influence of duplicates in the sample in comparison to deleting the duplicates. In this case, the depth value of the duplicated object should be higher based on the sample where the duplicates exists. The second property studies the impact of one single sample element on the resulting center-outward order. Here we consider two empirical probability measures. The first empirical probability measure is based on a sample that has an object which greatly differs from the other object. The second empirical probability measure is based on the same sample but without this greatly different object. An object $g_{diff}$ differs greatly from the other objects if this object has no attribute that any other object has. Or, equivalently, $g_{diff}$ and a further object of the sample $g_s$ are both elements of an extent $E$ if and only if $E = G$ is true.                                                                                                                                                                                                                       
Thus, this object $g_{diff}$ should have no impact on the center-outward order of the other objects. Note that there are scaling methods that exclude the existence of an object $g_{diff}$. In the Titanic example, see Table~\ref{tab:titanic_motivation}, we used the interordinal scaling method to include the numerical observation. This tells us that for two objects/passengers there must be at least one attribute that both objects have in common. Conversely, if ordinal scaling is used instead in the Titanic example, object $g_{diff}$ may occur.

\begin{enumerate}
	\item[(P9)] \textit{Respecting duplicates}: Let $\Kbb \in \varkappa$. Let $g_1, \ldots, g_n$ be a sample of $G$ with $n \in \mathbb{N}$. Assume that there exist $g_i$ and $g_j$ in the sample with $i\neq j$ which have identical attribute set (so $\Psi(g_i) = \Psi(g_j)$). We set 
	\begin{itemize}
		\item $\Pr^{(n)}$ to be the empirical probability measure of $g_1, \ldots, g_n$ with $\Pr^{(n)} \in \Prbb$, and
		\item $\Pr^{(n, -i)}$ to be the empirical probability measure $g_1, \ldots, g_{i-1},g_{i+1}, \ldots, g_{n}$ (without $g_i$) and $\Pr^{(n, -i)} \in \Prbb$. 
	\end{itemize}
	Then
	\begin{align*}
		D^{(n, -i)} (g_j, \Kbb) < D^{(n)}(g_j, \Kbb).
	\end{align*}
	\item[(P10)] \textit{Stability of the order}: Let $\Kbb \in \varkappa$ and let $g_1, \ldots, g_n$ be a sample of $G$ with  $n \in \mathbb{N}$. Assume that there exists an $i \in \{1, \ldots, n\}$ such that $g_i$ is greatly different from all other objects $g_1, \ldots, g_{i-1}, g_{i+1}, \ldots, g_n$ of the sample. This means that the only extents which contain $g_i$ and a further object $g_j$ with $j \in \{1, \ldots, i-1, i+1, \ldots, n\}$ is $G$. 
	Then the center-outward order of $g_1, \ldots, g_{i-1}, g_{i+1}, \ldots, g_n$ given by $D^{(n)}$ and $D^{(n, -i)}$ are the same. Recalling Definition~\ref{def:isomorph} and the notation in Property~(P9), we get
	\begin{align*}
		D^{(n)} (\cdot, \Kbb)_{\big| \{g_1, \ldots, g_{i-1}, g_{i+1}, g_n\}} \cong D^{(n, -i)} (\cdot, \Kbb)_{\big| \{g_1, \ldots, g_{i-1}, g_{i+1}, g_n\}}.
	\end{align*}
\end{enumerate} 

Let us take a look at Definition~\ref{def:allg_empir_depth_vorschrift} of the empirical depth function. There, the empirical probability measure is used for the definition. Thus, we start with a discussion on how duplicates or greatly different object effect the empirical probability measure of the extents. 
First, let us consider the case of Property~(P9) where two objects $g_1$ and $g_2$ are duplicates. This means that they have the same attributes and therefore lie in the same extent set. If the invariance on the attribute Property~(P2) is fulfilled, then $g_1$ and $g_2$ have the same depth. 
Furthermore, the probability reflects the duplicates for every extent which contains $g_1$ (and therefore also $g_2$).

Let us now consider Property (P10) and let $g_1, \ldots, g_n$ be a sample with $g_i = g_{diff}$ the greatly different object w.r.t.~the rest of the sample. One can observe that there exist three distinctive cases on how the empirical probability measures differ on $E \subsetneq G$ extent:
\begin{itemize}
	\item Case 1: $g_{diff}\in E$. Then $g_1, \ldots, g_i, g_{i+1}, \ldots, g_n \not\in E$ and  $\Pr^{(n)}(E) = 1/n \neq 0 = \Pr^{(n, -i)}(E)$.
	\item Case 2: $g_{\ell} \in E$ for some $\ell \in \{1,\ldots, i-1, i+1, \ldots, n\}$. Then $g_{diff} \not\in E$ and  $\Pr^{(n)}(E) = ((n-1)/n) \Pr^{(n, -i)}(E)$.
	\item Case 3:  $g_{\ell} \not\in E$ for all $\ell \in \{1, \ldots, n\}$. Then $\Pr^{(n)}(E) = 0 = \Pr^{(n, -i)}(E)$.
\end{itemize}
Thus, Cases 2 and 3 show that the order of the extents which do not contain $g_{diff}$ based on the probability measure $\Pr^{(n)}$ is the same as for $\Pr^{(n, -i)}$. Hence, the order structure given by the probability values of extents containing $g_1, \ldots, g_i, g_{i+1}, \ldots, g_n$ is not influenced by $g_{diff}$.

 These are observations on the empirical probability measure, but we only generally defined the (empirical) depth function. Thus, the question if a depth function fulfils these properties is concerned with how these connections between the empirical probability distributions are included in the precise definition of the mapping rule.

Properties (P9) and (P10) focus on two empirical probability measures and how their difference influences the empirical depth function. We end this subsection by discussing the behaviour of the depth function based on a sequence of empirical probability measures. Let $(\Pr^{(n)})_{n \in \mathbb{N}}$ be a sequence of empirical probability measures, with $\Pr^{(n)}$ being given by an independent and identical distributed (i.i.d.)~sample of $\Pr \in \Prbb$ with size $n \in \mathbb{N}$. By assumption we have that for all $n\in \mathbb{N}$, $\Pr^{(n)} \in \Prbb$. Property~(P11) discusses the consistency based on the empirical depth function towards the (population) depth function.
\begin{enumerate}
	\item[(P11)] \textit{Consistency}: Let $\Kbb \in \varkappa$ and $\Pr \in \Prbb$ be a probability measure on $G$. Let $\Pr^{(n)}$ be the empirical probability measure of an i.i.d.~sample $g_1, \ldots, g_n$ of $G$ with $n \in \mathbb{N}$ which is drawn based on $\Pr$. Then, 
	\begin{align*}
		\sup_{g \in G} \mid D^{(n)}(g, \Kbb) - D(g, \Kbb, \Pr)\mid \rightarrow 0 \text{ almost surely}.
	\end{align*}
\end{enumerate}

\subsection{Universality properties}

Finally, we introduce a notion of universality of depth functions. The main motivation is that in general there exists no strictly quasiconcave depth function, see~Theorem~\ref{prop:strength_strc_conc}. If one refrains from strict quasiconcavity, then a natural demand is to stick to quasiconcavity. However, the depth function that assigns every object the value zero is also quasiconcave but useless. 
The idea behind these properties is now the wish for a depth function to be \textit{as strictly quasiconcave as possible}. While the following motivation is based on the strictly quasiconcavity property, one can generalise this idea to all properties defined above. 

Before we introduce our notion of a richness of a depth function, we first indicate that defining this notion in a more naive way does
not lead to an adequate notion of richness. 
For now, let the formal context $\Kbb$ and probability measure $\Pr$ be fixed. An intuitive way to concretising the property of a quasiconcave depth function $D$ of being as strictly quasiconcave as possible is
to demand that there exists no other quasiconcave depth function $E$ that is more strictly quasiconcave than $D$. This can be formalised by saying that $E$ is more strictly quasiconcave than $D$ if everywhere, where $E$ violates strict quasiconcavity, also $D$ violates strict quasiconcavity. By violations of the strict quasiconcavity property, we mean that only equality and not strict inequality is fulfilled. More precisely, a depth function $E$ is more strictly quasiconcave than $D$ if
\begin{align*}
	quasiker(E(\cdot,\Kbb,\Pr))&\defeq \left\{ (A,g)\: \biggl| \: A\subseteq G, g \in G \text{ such that } g\in \gamma(A)\setminus A, \: E(g,\Kbb,\Pr) = \inf\limits_{\tilde{g}\in A} E(\tilde{g},\Kbb,\Pr)\right\} \\
	&\subsetneq quasiker(D(\cdot, \Kbb,\Pr)).
\end{align*}
In other words, there exists a pair $(A,g)\in quasiker(D(\cdot, \Kbb,\Pr))\setminus quasiker(E(\cdot, \Kbb,\Pr))$ with $g \in \gamma(A)$ such that the depth function $D$ violates for this pair the strictly quasiconcave property but $E$ does not. 

However, this definition leads to a problem.
The following situation can occur: Assume we have a quasiconcave depth function $D$ and $(A,g) \in quasiker(D(\cdot,\Kbb,\Pr))$ as well as $(\tilde{A},\tilde{g}) \in quasiker(D(\cdot,\Kbb,\Pr))$.
Then it may be the case that increasing the depth value $D(g,\Kbb,\Pr)$ by some small amount $\varepsilon$ leads to a depth function $E(\cdot,\Kbb,\Pr)$ that is still quasiconcave and that now fulfils $quasiker(E(\cdot,\Kbb,\Pr)) \subsetneq quasiker(D(\cdot,\Kbb,\Pr))$.
Now, assume that the same is true for increasing the depth of $\tilde{g}$ (without increasing the depth of $g$), but increasing the depth values of both $g$ and $\tilde{g}$ at the same time is not possible without violating the quasiconcavity property.
However, assume further that both the underlying probability measure $\Pr$ as well as the underlying context $\Kbb$ are perfectly symmetric w.r.t.~$g$ and $\tilde{g}$. In this case, it is not reasonable to arbitrarily increase one depth value to be more strictly quasiconcave. This means that according to the above definition of being as strictly quasiconcave as possible, there does not exist a possible and reasonable depth function that is as strictly quasiconcave as possible. This underlines the importance to emphasise that a depth function always inherently codifies both the structure of the underlying space as well as the concrete underlying probability measure.

We now introduce our notion of richness of a depth function. The following two universality properties are stated in a more general way such that they do not only capture the richness of  a depth function w.r.t.~the property of being quasiconcave. The introduced universality properties are also applicable to other properties of depth functions. Concretely, we take inspiration from a basic notion that is folklore in category theory. We adopt here the notion of an universal or a free object: Very roughly speaking, we say that a depth function $D$ is free w.r.t.~a set of some structural properties, if for every other depth function $E$ that obeys these properties, we can obtain $E(\cdot,\Kbb,\Pr)$ from $D$ as a composition of $D(\cdot,\Kbb,\tilde{\Pr})$ and an isotone function $f:\mathbb{R}\rightarrow \mathbb{R}$. Here, $\tilde{\Pr}$ is a suitable further probability measure. Informally, this means that we allow a meaningful depth function to have identical depth values if and only if the underlying probability measure and the structure of the space support this. 
 In other words, a depth function is free if it is flexible enough to imitate every other depth function (with the same properties) by supplying it with a suitable probability measure $\tilde{\Pr}$. With this idea in mind we now state two different notions of freeness, one weak and one strong notion. Therefore, let $Q\subseteq \{(P1),\ldots, (P10)\}$ be a set of properties from Sections~\ref{sec:representation_properties}, \ref{sec:order_preserving_properties} and~\ref{sec:empirical_properties} and recall Definition~\ref{def:isomorph}. A depth function $D$ is called

\begin{enumerate}
	\item[(P12)] \textit{weakly free} w.r.t.~a family $\mathscr{P}$ of probability measures on the object set $G$ and w.r.t.~$Q$, if it satisfies every property from $Q$ and if for every probability measure $\Pr$ and every depth function $E$  on $\Kbb$ that also satisfies all properties in $Q$, there exists a probability measure $\tilde{\Pr}\in \mathscr{P}$ and an isotone function $f:\mathbb{R}\longrightarrow \mathbb{R}$ such that 
	$$  f \circ D(\cdot, \Kbb, \tilde{\Pr}) \cong E(\cdot, \Kbb,\Pr).$$
	\item[(P13)] \textit{strongly free} w.r.t.~a family $\mathscr{P}$ of probability measures on the object set $G$ and w.r.t.~$Q$, if it satisfies all properties from $Q$ and if for every $\varepsilon >0$, there exists a class $\mathscr{P}^\varepsilon$ of probability measures from $\mathscr{P}$ with a diameter $d(\mathscr{P}^\varepsilon)$ lower than or equal to $\varepsilon$ such that for every other depth function $E$  on $\Kbb$ that also satisfies all properties in $Q$ and for every $\Pr \in \mathscr{P}$ there exists a probability measure $\tilde{\Pr} \in \mathscr{P}^\varepsilon$ and an isotone function $f:\mathbb{R}\longrightarrow \mathbb{R}$ with $$   f \circ D(\cdot, \Kbb, \tilde{\Pr}) = E(\cdot, \Kbb,\Pr).$$
	
	Here, the diameter of a set $\mathscr{P}$ of probability measures is defined as $d(\mathscr{P}) \defeq \sup_{P,Q \in \mathscr{P}}\sup_{A\in \Sigma} \lvert { P(A)-Q(A)} \rvert$ where $\Sigma$ is the assumed underlying $\sigma$-field.
\end{enumerate}

\begin{remark}
	We say that a depth function $D$ is richer than a depth function $E$ if and only if we can modify $D$ in such a manner that the center-outward order given by $E$ can be reconstructed by $D$. This means that we can find a probability measure $\tilde{\Pr}$ and an isotone function $f$ such that $f \circ D(\cdot, \Kbb, \tilde{\Pr})$ equals $E$. Note that $f$ only needs to be isotone and not bijective. 
	 Therefore, $D$ must be at least as rich as $E$ in distinguishing the individual objects of $G$. Finally, observe that $\tilde{\Pr}$ and $f$ depend on $Q, \Pr, \Kbb, D$ and $E$.
\end{remark}

Within the notion weakly free we completely detach from $\Pr$ by allowing $\tilde{\Pr}$ to be any probability measure from $\mathscr{P}$.
For the notion strongly free we restrict the allowed probability measures $\tilde{\Pr}$. Therefore, strongly free implies weakly free.
The idea behind restricting the probability measure $\tilde{\Pr}$ is the requirement that a depth function should be able to detach itself to some degree from the underlying probability measure $\Pr$. More precisely, there exists a balanced measure $\Pr^*$ such that any structure within the depth values that is not due to the structure of the formal context can already be modified by an arbitrarily small modification of this balanced measure $\Pr^*$. The reason why we did not explicitly use such a balanced measure $\Pr^*$ in the definition is a technical one. The balanced measure does not need to be a probability measure (e.g., there is no uniform distribution on the real numbers). Therefore we work here with the set $\mathscr{P}^\varepsilon$ with arbitrary $\varepsilon > 0$.

\begin{theorem}\label{strongly_free_exists}
	For the hierarchical nominal context given in Table~\ref{crosstable_hierarchical_nominal} of Example \ref{hierarchical nominal}, there exists a depth function $D$ that is strongly free w.r.t.~quasiconcavity (and w.r.t.~the family of all probability measures on $G$), namely the function $D$ given by:

	\begin{align*}
		D(g,\Kbb,\Pr)&=\begin{cases}
			1 & \mbox{ if } g=\;g^* \defeq\argmax\limits_{\tilde{g} \in G}\Pr(\{\tilde{g}\})\\
			1/2 & \mbox{if on the first level, $g$ has the same category as $g^*$}\\
			0& \mbox{else }
		\end{cases}.
	\end{align*}
	Here we assume that the argument of the maximum is unique. Otherwise we arbitrarily choose one of the arguments of the maximum.
	
\end{theorem}
\begin{proof}
	First observe that $D$ is indeed quasiconcave, because the contour sets are extent sets. The reason is that the images of $\gamma$ are exactly the one element sets, the whole set $G$ and the two element sets $\{ a_1a_2,a_1b_2\}$ and $\{b_1a_1,b_1a_2\}$. Now, for $\varepsilon >0$ define the class of probability measures $\mathscr{P}^\varepsilon \defeq \left\{\Pr \mid \forall g \in G \text{ we have } \Pr(\{g\}) \in [1/4 -\varepsilon/4, 1/4 +\varepsilon /4 ]\right\}.$
	This class has in fact diameter lower than or equal to $\varepsilon$. Now, let $\varepsilon >0 $, let $\Pr$ be an arbitrary probability measure on $G$, and let $E$ be an arbitrary quasiconcave depth function. Let further $g$ be the object with the highest depth according to $E$. Without loss of generality we assume that $g =a_1a_2$. Additionally, we assume that there exists exactly one $g$ with maximal depth (otherwise, arbitrarily choose from the objects with maximal depth). Then define $\tilde{\Pr}\in \mathscr{P}^\varepsilon$
	by $\tilde{Pr}(\{g\}) \defeq 1/4+\varepsilon$ and $\tilde{\Pr}(\{\tilde{g}\}) \defeq 1/4-\varepsilon/3$ for $\tilde{g}\neq g$. If $\varepsilon >3/4$, we set $\varepsilon $ to $3/4$. Because $E$ is assumed to be quasiconcave, and since the contour sets are always nested, the contour sets of $E$ are exactly the sets $\{a_1a_2\}, \{a_1a_2,a_1b_2\}$ and $G$. (For $E$ it could also be the case that the set of contour sets is only a subset of these three sets). 
	At the same time, due to construction, the contour sets of $D$ are exactly the same which means that we can construct an isotone function $f:\mathbb{R}\rightarrow \mathbb{R}$ with $f\circ D(\cdot, \Kbb,\tilde{\Pr})=E(\cdot,\Kbb,\Pr)$. Since in particular $\varepsilon>0$ was arbitrary, $D$ is in fact strongly free w.r.t.~quasiconcavity and w.r.t.~the family of all probability measures on $G$.
	
\end{proof}

	Also for a finite hierarchical nominal context with more than two levels and with more than two categories in every level, one can show that there exists a depth function that is strongly free w.r.t.~quasiconcavity and w.r.t.~the family of all probability measures.

\begin{remark}
	Note that while for example quasiconcavity implies starshapedness, freeness w.r.t.~quasiconcavity generally does not imply freenness w.r.t.~starshapedness. This is because a quasiconcave depth function $D$ is not able to imitate a depth function $E$ that is starshaped but not quasiconcave, if such a depth function $E$ exists at all. In our concrete Example~\ref{hierarchical nominal} of the hierarchical nominal context with two levels and two categories on every level, quasiconcavity is accidentally equivalent to starshapedness. For more than two levels this is not the case anymore.
\end{remark}

\section{Example: Generalised Tukey depth function \label{sec:tukeys_depth}}
In Section~\ref{sec:des_properties}, we generally introduced structural properties for depth functions based on FCA.  
One purpose of these structural properties is to give a systematic basis to analyse depth functions using FCA. Therefore, in this section, we take the opportunity to analyse the generalised Tukey depth given in Section~\ref{sec:definition}. Before we begin with discussing the structural properties, we look at some general aspects of the generalised Tukey depth. Here and in the following, we assume that $G$ is a set, $\varkappa$ a set of formal contexts and $\Prbb$ a set of probability measures on $\sigma$-fields which contain every extent set given by $\varkappa$.

In Definition~\ref{def:gen_tukey} the generalised Tukey depth is defined by using the extent sets of single attributes. Thus, in the computation one apparently only takes a proper subset of all possible extents into account, see, e.g., Example~\ref{examp:tukey_rechenbsp} and~\ref{examp:posets_tukey_depth}. The reason for this is that considering all extent sets instead of only those induced by a single attribute does not change the depth value. Compare to Section~\ref{sec:definition} where this has already been mentioned.

\begin{theorem}\label{theo:tukey_extent}
	Let $G$ be an object set. For every formal context $\Kbb$ with object set $G$ and every probability measure $\Pr$ on $G$ we get for the generalised Tukey depth function
	\begin{align*}
		T(g, \mathbb{K}, \Pr) = 1 - \sup_{m \in M \setminus \Psi(g)}\Pr(\Phi(m)) = 1 - \sup_{B \subseteq M \setminus \Psi(g)}\Pr(\Phi(B)).
	\end{align*}
\end{theorem}
\begin{proof}
	We have to show that 
	\begin{align}\label{eq:proof_tukey_extent}
	\sup_{m \in M \setminus \Psi(g)}\Pr(\Phi(m)) = \sup_{B \subseteq M \setminus \Psi(g)}\Pr(\Phi(B)).
	\end{align}
	Since $B \subseteq M \setminus \Psi(g)$ is a superset of $m \in M \setminus \Psi(g)$ we immediately get $\le$ in Equation~(\ref{eq:proof_tukey_extent}). For the reverse inequality let $B \subseteq M \setminus \Psi(g)$ be arbitrary. Then for every $m \in B$ we get that $\Phi(m) \supseteq \Phi(B)$ and thus $\Pr(\Phi(m)) \ge \Pr(\Phi(B))$. Note that $m \in M \setminus \Psi(g)$. Together with the properties of a supremum, we get $\ge$ in Equation~(\ref{eq:proof_tukey_extent}).
\end{proof}

\begin{remark}\label{rem:tukey_extent}
	Analogously to Theorem~\ref{theo:tukey_extent} one can show that $$T(g, \mathbb{K}, \Pr) = 1 - \sup_{m \in M \setminus \Psi(g)}\Pr(\Phi(m)) = 1 - \sup_{\substack{A \subseteq G\setminus\{g\}\\ A \text{ extent} }}\Pr(A)$$
	is true for $g \in G$, $\Kbb$ a formal context and probability measure $\Pr$. Again, we obtain $\le$ immediately. For the reverse, assume for simplicity that the supremum is attained at extent $A \subseteq G\setminus \{g\}$. Since $g \not\in A$ holds, there exists $m \in \Psi(A)$ such that $g \not\in \Phi(m)$. In particular, we get for this $m$ that $\Pr(\Phi(m))\ge \Pr(A)$ is true which gives indeed equality above.
\end{remark} 
\begin{remark}\label{rem:tukey_margin}
Theorem~\ref{theo:tukey_extent} shows that the generalised Tukey depth uses only the marginal distribution over the attributes. Thus, it is possible to change the incidence relation $I$ and thereby change the dependency structures of the objects, but the generalised Tukey depth remains the same. Obviously, this statement is particularly true if the marginal distribution remains the same. However, since only the supremum of the extents sets of single attributes is used, even the marginal distribution can change and still the generalised Tukey depth remains the same. An example is given in Table~\ref{tab:limits_tukey_allP8} (middle and right). Here, we have $T(g, \Kbb_1, \Pr) = T(g, \Kbb_2, \Pr)$ for all $g \in \{g_1, g_2, g_3\}$ (with uniform probability measure on $\{g_1, g_2, g_3\}$), although in the second formal context, we have that $g_2$ implies $g_3$.
\end{remark}	
	\begin{table}[!htb]
	\begin{minipage}{.3\linewidth}
	\centering
	\begin{tabular}{lccc}
        		\hline
 			& $m_1$ & $m_2$ & $m_3$\\\hline
			$g_1$ & $\times$ &   &   \\
			$g_2$ &  &$\times$  & \\
			$g_3$ & & & $\times$\\ 
			\hline
     \end{tabular}
	\end{minipage}
    \begin{minipage}{.3\linewidth}
      \centering
	\begin{tabular}{lcccc}
        		\hline
 			& $m_1$ & $m_2$ & $m_3$ & $m_4$\\\hline
			$g_1$ & $\times$ &   &  & $\times$  \\
			$g_2$ & & $\times$ &$\times$  &  \\
			$g_3$ & & & $\times$& $\times$\\ 
			\hline
     \end{tabular}
    \end{minipage}%
    \begin{minipage}{.3\linewidth}
      \centering
        \centering
	\begin{tabular}{lcccc}
        		\hline
 			& $m_1$ & $m_2$ & $m_3$& $m_4$\\\hline
			$g_1$ & $\times$ &   & & $\times$  \\
			$g_2$ & & $\times$ &$\times$  & \\
			$g_3$ & & $\times$ & $\times$ & $\times$\\ 
			\hline
     \end{tabular}
    \end{minipage} 
    \caption{For the formal context on the left Property~(P8) is never true (see Section~\ref{sec:order_preserving_properties}). The middle ($\Kbb_1$) and right ($\Kbb_2$) formal contexts show the effect of the supremum in the definition of the generalised Tukey depth.}
    \label{tab:limits_tukey_allP8}
\end{table}


\subsection{Representation properties}
The first two structural properties, the representation properties, aim to ensure that the structure in the extent set is reflected in the depth function. The generalised Tukey depth function is based on the set of extents, see Theorem~\ref{theo:tukey_extent} and Remark~\ref{rem:tukey_extent}. More precisely, for each $g$ the depth is based on the probabilities of those extents which do not contain~$g$. Since the function $i$ in Property (P1) is bijective and preserves extents and probabilities, this implies that the depth values are preserved as well. This shows that Property (P1) holds. Furthermore, since two objects which equal in their attributes, always lie in the same extent sets, they have to have the same depth values. Thus, Property (P2) is true for the generalised Tukey depth.

\subsection{Order preserving properties}

In Section~\ref{sec:des_properties}, the order preserving properties are ordered in their strength. Thus, we prove that the contourclosed property (P7i) holds and with this further order preserving properties follow by Theorem~\ref{theo:P7i_ii_equi_P8_impl_P7} and \ref{lem:P5_implies_P3_4}.

\begin{theorem}
	The generalised Tukey depth fulfils Property (P7i).
\end{theorem}
\begin{proof}
	We prove Property (P7i) by contradiction. Assume that there exists a formal context $\Kbb$ on an object set $G$ and a probability measure $\Pr$ together with an $\alpha \in im(T(\cdot, \Kbb, \Pr))$ such that $Cont_{\alpha}$ is not an extent set. This means that $\gamma(Cont_{\alpha}) \supsetneq Cont_{\alpha}$ since $\gamma$ is a closure operator. Thus, there exists $g \in \gamma(Cont_{\alpha})\setminus Cont_{\alpha}$. The generalised Tukey depth, see Definition~\ref{def:gen_tukey}, is based on the extent sets induced by a single attribute that the object $g$ does not have. Since $g \in \gamma(Cont_{\alpha})$ we know that for every attribute which $g$ does not have there exists at least one $\tilde{g} \in Cont_{\alpha}$ which does not have this attribute either. Hence, we have for all $ m \in M\setminus \Psi(g)$ there exists $\tilde{g}\in Cont_{\alpha}$ such that $1- \Pr(\Phi(m)) \ge 1- \sup_{\tilde{m} \in M \setminus \Psi(\tilde{g})}\Pr(\Phi(m)).$
By definition of the contour set $\alpha$, we get $T(g, \Kbb, \Pr) = 1 - \sup_{m \in M \setminus \Psi(g)}\Pr(\Phi(m)) \ge \inf_{\tilde{g} \in Cont_{\alpha}}\left(1- \sup_{\tilde{m} \in M \setminus \Psi(\tilde{g})}\Pr(\Phi(\tilde{m}))\right) = \inf_{\tilde{g} \in Cont_{\alpha}}T(\tilde{g}, \Kbb, \Pr).$
This is contradicting $g \not\in Cont_{\alpha}$ and we can conclude that Property (P7i) holds for every formal context $\Kbb$ and probability measure $\Pr$.
\end{proof}
	
\begin{remark}
	Now we apply Theorem~\ref{lem:P5_implies_P3_4} and \ref{theo:P7i_ii_equi_P8_impl_P7}, see overview in Figure~\ref{fig:uebersicht_property_structure}, and obtain that since Property (P7i) is true for the generalise Tukey depth, the following properties hold as well: minimality property (P3), maximality property (P4), isotonicity property (P5), starshapedness property (P6) and the quasiconcavity property (P7ii).
\end{remark}

In Section~\ref{sec:des_properties} we showed that there exist formal contexts such that the strictly quasiconcavity property (P8) does not hold for any depth function and any probability measure. Hence, here we are interested in formal contexts such that the generalised Tukey depth is strictly quasiconcave. A small example of a formal context such that the quasiconcavity is fulfilled is given in Table~\ref{tab:fc_tukey_p8_notP10_fulfilled} (left).

\begin{table}[!htb]
    \begin{minipage}{.5\linewidth}
      \centering
	\begin{tabular}{lcccc}
        		\hline
 			& $m_1$ & $m_2$ & $m_3$ & $m_4$\\\hline
			$g_1$ & $\times$ & $\times$ &$\times$ &\:   \\
			$g_2$ & $\times$ &$\times$  & \: & $\times$\\
			$g_3$ & $\times$& $\times$ & $\times$ & $\times$\\ 
			\hline
     \end{tabular}
    \end{minipage}%
    \begin{minipage}{.5\linewidth}
      \centering
        \centering
	\begin{tabular}{lcccc}
        		\hline
 			& $m_1$ & $m_2$ & $m_3$ \\\hline
			$g_1$ & $\times$ & \: & \:  \\
			$g_2$ & \:  &$\times$  & \: \\
			$g_3$ & \: & $\times$ & $\times$ \\ 
			\hline
     \end{tabular}
    \end{minipage} 
    \caption{For the formal context on the left holds Property (P8) for the generalised Tukey depth. The formal context does not fulfil Property (P10) for the generalised Tukey depth.}
    \label{tab:fc_tukey_p8_notP10_fulfilled}
\end{table}


To check this claim, we take a look at Property (P8) applied on the generalised Tukey depth. Thus, Property~(P8) is satisfied if and only if for all $A \subseteq G$ we have for all $\tilde{g} \in \gamma(A)\setminus A$
\begin{align}\label{ineq: tukey_p8_fulfilled}
	1- \sup_{m \in M \setminus \Psi(\tilde{g})}\Pr(\Phi(m)) > \inf_{g \in A} \left( 1- \sup_{m \in M \setminus \Psi(g)}\Pr(\Phi(m))\right)
	\Leftrightarrow  \sup_{m \in M \setminus \Psi(\tilde{g})}\Pr(\Phi(m)) < \sup_{g \in A} \sup_{m \in M \setminus \Psi(g)}\Pr(\Phi(m)).
\end{align}
Note that $\le$ in Inequality~(\ref{ineq: tukey_p8_fulfilled}) follows by Property~(P7ii).
To check whether the generalised Tukey depth based on the formal context given in Table~\ref{tab:fc_tukey_p8_notP10_fulfilled} (left) is strictly quasiconcave, we go through every subset $A \subseteq G$  with $\gamma(A)\setminus A\neq \emptyset$. This is only true for $A = \{g_1, g_2\}$ with $\gamma(A) = \{g_1, g_2, g_3\}$. The attributes of $g_3$ are a superset of the union of the attributes of $g_1$ and $g_2$ (i.e. $\Psi(\{g_3\}) \supseteq \Psi(\{g_1\}) \cup \Psi(\{g_2\})$). Observe that the attributes of $g_1$ and $g_2$ which maximise the supremum of the corresponding generalised Tukey depth function are attributes of $g_3$. Thus, when ensuring that there exists an attribute $m \in \Psi(g_3)\setminus \Psi(g_2) \cap \Psi(g_1)$ such that $\Pr(\Phi(m))$ is strictly larger than $\Pr(\Phi(\tilde{m}))$ based on those attributes $\tilde{m}$ which are not an element of $\Psi(g_3)$, the generalised Tukey depth is strictly quasiconcave. Note that it is sufficient that the depth of either $g_1$ or $g_2$ must be below the depth of $g_3$. Furthermore, we want to point out that we also have to assume that there does not exist a further object $g_4$ which is a duplicate of $g_3$, because then $\Psi(g_3)\setminus \Psi(g_2) \cap \Psi(g_1) = \emptyset$. With this, we get the strict part in Inequality~(\ref{ineq: tukey_p8_fulfilled}).
Based on the idea above, we can say that the generalised Tukey depth satisfies Property (P8) for every formal context which is an element of the following set:
 
\begin{align*}
  \Cscr^{P8} = \left\{ \begin{array}{l} (G,M,I) \text{ formal context without}\\ \text{duplicates according to attributes}  \end{array} \bigg| \begin{array}{l} \text{for every subset } A \text{ and } g \in \gamma(A)\setminus A \text{ we have } 
  \Psi(g) \supseteq \cup_{a \in A} \Psi(a) \text{ and }\\
  \exists m \in \Psi(g)\setminus \cap_{a \in A} \Psi(a)\colon \Pr(\Phi(m)) > \sup_{m \in M \setminus \Psi(g)} \Pr(\Phi(m))
  \end{array} \right\}.
\end{align*}



\subsection{Empirical (sequence) properties}
The discussion about empirical (sequence) properties profits from the strong impact of the (empirical) probability measure on the definition of the generalised Tukey depth, see Definition~\ref{def:gen_tukey} and~\ref{def:empirical_tukey}. With this we can immediately follow that Property~(P9) is fulfilled by the generalised Tukey depth. The stability of the order property does not hold. Consider again the formal context given by Table~\ref{tab:fc_tukey_p8_notP10_fulfilled} (right). Then $g_1$ is an outlier w.r.t.~$\{g_2, g_3\}$. 
For the empirical generalised Tukey depth based on $\Pr^{(2)}$ for $\{g_2, g_3\}$ we get $T^{(2)}(g_3, \Kbb) =1 > 1/2 = T^{(2)}(g_2, \Kbb)$. When adding object $g_1$ to the probability measure $\Pr^{(3)}$ we get that  $T^{(3)}(g_3, \Kbb)= 2/3 = T^{(3)}(g_2, \Kbb)$. The order of $g_2$ and $g_3$ is not stable. Since this can always occur for small $n \in \mathbb{N}$, we cannot restrict the set of probability measures or formal contexts.

It remains to discuss the consistency property (P11). One can show that if the set of extents induced by single attributes $\mathcal{ES} \defeq \{\Phi(m)\mid m \in M\}$
is a Glivenko-Cantelli class w.r.t.~the underlying probability measure $\Pr$, see \cite{dudley91}, then the generalised Tukey depth is indeed consistent. This, for example, is given when $\mathcal{ES}$ has a finite VC dimension. 
\begin{theorem}\label{consistency}
	Let $\Kbb$ be a formal context with $G$ as object set. Let $\text{Pr}^{(n)}$ be the empirical probability measure given by an i.i.d.~sample of sizes $n \in \mathbb{N}$ and based on probability measure $\Pr$. When $\mathcal{ES}$ is a Glivenko-Cantelli class w.r.t.~$\Pr$, then $T^{(n)}$ is consistent and Property (P11) is fulfilled with
	\begin{align*}
		\sup_{g \in G}\: \left| T^{(n)}(g, \Kbb)- T(g, \Kbb, \Pr) \right|\: \overset{n \to \infty}{\longrightarrow} \:0  \: \text{ almost surely}.
	\end{align*} 
\end{theorem}
\begin{proof}
	The proof is similar to the proof given by \citep[p. 1816f]{donoho92}. Since we assume that $\mathcal{ES}$ is a Glivenko-Cantelli class w.r.t.~$\Pr$, we get $\sup_{E \in \mathcal{ES}} \left| \text{Pr}^{(n)}(E) - \Pr(E) \right| \:  \overset{n \to \infty}{\longrightarrow} \:0$ almost surely.
	With this, the claim follows from
	\begin{align*}
		&\sup_{g \in G}\: \left| 1- \sup_{m \in M\setminus\Psi(g)} \text{Pr}^{(n)}(\Phi(m)) - \left(1- \sup_{m \in M\setminus\Psi(g)} \Pr(\Phi(m))\right)  \right| 
		= \sup_{g \in G}\: \left| \sup_{m \in M\setminus\Psi(g)} \text{Pr}^{(n)}(\Phi(m)) - \sup_{m \in M\setminus\Psi(g)} \Pr(\Phi(m)) \right|\\
		\le &\sup_{E \in \mathcal{ES}} \left| \text{Pr}^{(n)}(E) - \Pr(E) \right| \:  \overset{n \to \infty}{\longrightarrow} \:0  \: \text{ almost surely.}
	\end{align*}
\end{proof}
\begin{remark}
	If $\mathcal{ES}$ has a finite VC dimension, then the convergence stated in Theorem~\ref{consistency} is additionally uniform over all possible underlying probability measures. A finite VC dimension of $\mathcal{ES}$ is given, for example, if $M$ is finite.
\end{remark}

\begin{remark}
	We know that the generalised Tukey depth fulfils Property (P1). Let us assume that for a formal context $\tilde{\Kbb}$ the set of all extents induced by one single attribute is not a Glivenko-Cantelli class w.r.t.~some probability measure $\tilde{\Pr}$. Then, in some cases, one can define a second formal context ${\Kbb}$ and probability measure $\Pr$ such that there exists a bijective and bimeasurable function $i$ which preserves the extents and the probability measure. If the extents based on a single attribute given by ${\Kbb}$ is a Glivenko-Cantelli class w.r.t.~probability measure $\Pr$, then we can transfer the consistency based on ${\Kbb}$ onto $\tilde{\Kbb}$. An example for this is given by the following two formal contexts: Let $G= \mathbb{R}^d$ and ${\Kbb}$ be the formal context defined in Example~\ref{exampl: spatial_1}. For $\tilde{\Kbb}$ let $\tilde{M}$ be the set of all topologically closed convex sets and $\tilde{I}$ the binary relation with $(g,\tilde{m}) \in I$ if and only if $g$ is an element of the corresponding convex set $\tilde{m}$. The extents equal again all topologically closed convex sets. Thus by applying Property (P1) we can replace $\tilde{\Kbb}$ with $\Kbb$ and can prove that also for $\tilde{\Kbb}$ the generalised Tukey depth is consistent. 
\end{remark}

\subsection{Universality properties}
Let us end this discussion on the generalised Tukey depth by considering the universality properties. The next two theorems prove that the generalised Tukey depth is weakly free but not strongly free w.r.t.~the quasiconcavity property. Therefore, recall the hierarchical nominal example, see Example~\ref{hierarchical nominal} and~\ref{examp:hierarchical_nominal_tukey}.

\begin{theorem}
	Let $G$ be finite. Then the generalised Tukey depth is weakly free w.r.t.~the family of all probability measures on $G$ and w.r.t.~the quasiconcavity property (P7ii).
\end{theorem}

\begin{proof}
	Let $G$ be finite, let $\Kbb$ be an arbitrary context with object set $G$ and let $\Pr$ be an arbitrary probability measure on $G$. Moreover, let $E$ be an arbitrary quasiconcave depth function. We set $E_1,\ldots, E_K$ to be the unique increasingly ordered depth values of $E(\cdot,\Kbb,\Pr)$ and let $k_i$ be the number of objects with depth value $E_i$. Furthermore, let $G_i$ denote the boundaries of the contour sets, i.e., the set of objects with depth $E_i$. We now have to construct a probability measure $\tilde{\Pr}$ and an isotone function $f$ such that $f \circ T(\cdot,\Kbb,\tilde{\Pr})= E(\cdot,\Kbb,\Pr)$. 
	
	To achieve this, we first consider properties this probability function needs to fulfil. Let $\tilde{\Pr} \in \Prbb$ be a possible candidate. We set	
	$p_i\defeq\tilde{\Pr}(G_i)$ and $p_i^{min}\defeq \min_{g \in G_i}\tilde{\Pr}(\{g\})$. Since $E$ is quasiconcave, we have $g \notin \gamma (G_{\ell+1}\cup \ldots \cup G_K)$ for all $g \in G_{\ell}$ with $\ell \in \{1,\ldots,K-1\}$. This implies that for every $\ell \in \{1,\ldots,K-1\}$ and every $g \in G_{\ell}$ there exists an attribute $m\in M$ with $(g,m)\notin I$, but $(h,m) \in I$ for all $h \in G_{\ell+1}\cup \ldots \cup G_{K}$. Therefore, we get for $g \in G_{\ell}$
	\begin{align}\label{inequality_layers}
		p_{\ell}^{min} &\leq T(g,\Kbb,\tilde{\Pr}) \leq \sum_{i=1}^{\ell} p_i.
	\end{align}
	Now, we construct probability measure $\tilde{\Pr}$ needed for Property (P12) in three steps: First set $\tilde{\Pr}(g) =1$ for $g \in G_1$. Then recursively set $\tilde{\Pr}(g) = \sum_{i=1}^{\ell} \tilde{\Pr}(G_i) +1$ for $g \in G_{\ell+1}$. Finally, we normalise $\tilde{\Pr}$ by $\tilde{\Pr}(g) \defeq \tilde{Pr}(g)/\sum_{g \in G}\tilde{\Pr}({g})$ to get a probability measure. Due to construction and Inequality~(\ref{inequality_layers}), it is ensured that $T(g,\Kbb,\tilde{Pr}) < T(h ,\Kbb,\tilde{Pr})$ for $g \in G_{\ell}$ and $h \in G_{\ell+1}$. 
	Thus, one can define function $f$ by setting $f(T(g,\Kbb,\tilde{\Pr}))=E_{\ell}$ for $g \in G_{\ell}$ and isotonically extending it to function $f:\mathbb{R}\rightarrow\mathbb{R}$. With this, we get
	$f \circ T(\cdot,\Kbb,\tilde{Pr})=E(\cdot ,\Kbb,\Pr)$ which shows that claim. 
\end{proof}

\begin{theorem}\label{theo: tukey_strongly_free}
	For the context describing hierarchical nominal data given in Table~\ref{crosstable_hierarchical_nominal} the generalised Tukey depth is not strongly free w.r.t quasiconcavity and w.r.t.~the family of all probability measures on $G$.
\end{theorem}

\begin{proof}
	First note that the depth function $D$ from Theorem \ref{strongly_free_exists} is quasiconcave. Furthermore, $D$ is flexible enough to assign every arbitrary object $g$ the highest depth $1$. This implies that the object $h$ with the same category on Level $1$ as $g$ has the second highest depth $1/2$ and all other objects have depth $0$. Now, the generalised Tukey depth $T(\cdot,\Kbb,\Pr)$ needs to be as flexible as $D$ to be strongly free w.r.t.~quasiconcavity.
	Looking at Table~\ref{crosstable_hierarchical_nominal}, first observe that for object $g=a_1a_2$ to have the highest generalised Tukey depth, it is necessary that $\Pr(\{a_1a_2,a_1,b_2\})\geq 0.5$. This is because if we have $\Pr(\{a_1a_2,a_1,b_2\}) <  0.5$, we can conclude that $T(a_1a_2,\Kbb,\Pr) \leq 1- \Pr(\{b_1a_2,b_1b_2\}) =\Pr(\{a_1a_2,a_1,b_2\})<0.5$. Since either $\Pr(\{b_1a_2\})$ or $\Pr(\{b_1b_2\})$ is smaller or equal to $ 0.5$, we get that one of the objects $b_1a_2$ or $b_1b_2$ has a generalised Tukey depth of larger or equal to $0.5$. Now, for object $h=a_1b_2$ to have the second highest depth, it is necessary that the corresponding supremum in Definition~\ref{def:gen_tukey} is attained for an attribute that differs from $g=a_1a_2$. Because both $g$ and $h$ do not have attribute $b_1$, it is clear that attribute $b_1$ is relevant for the supremum within the generalised Tukey depth both for objects $g$ and $h$. For $h$ to have the second highest depth, it is therefore necessary that for $h$ the supremum is attained for attribute $a_1a_2$. But for this it is necessary that 	
	\begin{align}
		\Pr(\{a_1b_2\})\geq \Pr(\{b_1a_2,b_1b_2\}).\label{ineq_freenessproof}
	\end{align}	
	Now, let $\varepsilon=0.1$ and assume that  $T$ is strongly free w.r.t.~quasiconcavity. Thus, we can find a family $\mathscr{P}^\varepsilon$ of diameter $\varepsilon=0.1$ with the corresponding properties. Now take $\Pr^*\in \mathscr{P}^\varepsilon$ such that the generalised Tukey depth has the highest depth at $g=a_1a_2$ and for $h=a_1b_2$ the second highest depth. With the above considerations we can conclude first that $\tilde{Pr}(\{a_1a_2,a_1b_2\})\geq 0.5$ and, thus, $\Pr^*(\{a_1a_2,a_1b_2\}) \geq 0.5 -\varepsilon =0.4.$ Additionally, because also object $b_1a_2$ could be the object with highest depth w.r.t.~$D$ and some probability measure $\Pr$, an analogous argumentation implies that also $\Pr^*(\{b_1a_2,b_1b_2\}) \geq 0.5 -\varepsilon =0.4$. But this, together with Inequality (\ref{ineq_freenessproof}) implies that $\Pr^*(a_1b_2)\geq 0.4$. With analogous argumentations, one can also show $\Pr^*(\{g\})\geq 0.4$ for all other objects. But this is a contradiction to $\Pr^*$ being a probability measure because $\Pr^*(G)=0.4*4=1.6$. This provides the claim.
	
\end{proof}

\section{Conclusions}\label{sec:conclusion}

This article introduced a general notion of depth functions for non-standard data based on FCA. This covers the analysis of centrality and outlyingness for a wide variety of different data types. In order to enable a discussion, especially about these non-standard data, which are not given in a standard statistical data format, we used FCA and introduced structural properties. 
In addition to adopting the properties discussed in $\mathbb{R}^d$, we also addressed the issue that non-standard data may inherit a central-outward order that should also be reflected in the depth function. We also provided a framework for analysing depth functions based on FCA. Finally, we used this framework to discuss and analyse the generalised Tukey depth.
Building on this, there are several promising areas for further research, including (but not limited to)

\textbf{Further concrete mapping rules for depth functions:} In Theorem~\ref{strongly_free_exists} we defined a depth function $D$ which is strongly free w.r.t.~Property (P8) and the formal context defined in Table~\ref{crosstable_hierarchical_nominal}. One can show that this depth function $D$ only uses the probabilities of Level 1 to obtain the object $g$ with the highest depth value. In further research, it is of interest to define mapping rules for depth functions which take more levels into account. We already introduced the union-free generic depth function for the special case of $G$ being the set of partial orders, see \cite{blocher23}. This is an adaptation of the simplicial depth in $\mathbb{R}^d$, see~\cite{liu90}, to the set of all partial orders. As shortly denoted there this can be generalised to general formal contexts. \cite{hu23} give another examples for a depth function using FCA. The authors there propose concrete outlier measures based on formal concept analysis. Since there is a strong relationship between outlier measures and depth functions, the properties here can also be used to analyse this measure.

\textbf{Larger scale application:} In this article the focus lies on the general analysis of depth functions and the generalised Tukey depth. Hence, the examples discussed here aim to support and clarify claims of theoretical/structural properties of depth functions. In further research, a discussion based on the perspective of a larger scale application is of interest. For example, in \citep{blocher23} we applied the union-free generic depth function, which can be embedded into our concept of depth upon partial orders that represent the performance of machine learning algorithms.

\textbf{Statistical inference:} Building on this not only a descriptive analysis, but also statistical inference methods can be developed. Since depth functions describe in a robust and nonparametric manner the order of the data this can be used to define statistical tests and models. Analogously to the approach here, where we, among others, transferred already existing properties to our concept, this can be done for statistical inference methods, see~\cite{li04} as starting point.

\textbf{Analysis of generalised Tukey depth based on one specific scaling method:} Here, we generally analysed the generalised Tukey depth function without focus on one specific scaling method. In \cite{blocher22} we discussed the special case of $G$ being the set of partial orders and one scaling method. Since one fixed scaling method on a set $G$ gives us more structure on the corresponding closure system an analysis can lead to more structural properties of the depth function.

\section*{Acknowledgments}
HB and GS sincerely thank Thomas Augustin for the many helpful suggestions during the preparation of this article. They are also grateful to the two reviewers and the associated editor for their very important suggestions, which, among others, helped to improve the accessibility of the article. HB and GS gratefully acknowledge the financial and general support of the LMU Mentoring program. HB sincerely thanks Evangelisches Studienwerk Villigst e.V. for funding and supporting her doctoral studies. 


%
%
%

\bibliographystyle{myjmva}
\bibliography{trial.bib}
\end{document}